\input amstex
\input epsf
 \documentstyle{amsppt}
\loadeurm

\font\eurs=eurm7

\def\S{S}
\def\g{\frak g}
\def\h{\frak h}
\def\n{\frak n}
\def\b{\frak b}

\def\fy{\varphi}
\def\R{\Bbb R}
\def\F{{\Bbb F}}
\def\Z{\Bbb Z}
\def\C{\Bbb C}

\def\ZZ{\Cal Z}
\def\AA{\Cal A}
\def\X{\Cal X}
\def\L{\Cal L}
\def\P{\Cal P}
\def\D{{\widehat{D}}}
\def\ls{\leqslant}
\def\gs{\geqslant}
\def\one{\mathbf 1}

\def\dim{\operatorname{dim}}
\def\diag{\operatorname{diag}}
\def\deg{\operatorname{deg}}
\def\rank{\operatorname{rank}}
\def\corank{\operatorname{corank}}
\def\Spec{\operatorname{Spec}}
\def\Ad{\operatorname{Ad}}
\def\Norm{\operatorname{Norm}}
\def\End{\operatorname{End}}
\def\Trace{\operatorname{Tr}}
\def\Mat{\operatorname{Mat}}
\def\s{\operatorname{sign}}
\def\one{{\bold 1}}
\def\el{\text{\eurs l}}
\def\er{\text{\eurs r}}
\def\gr{\Cal E}

\def\ii{\bold i}
\def\jj{\bold j}
\def\kk{\bold k}
\def\cC{{\AA(Z)}}
\def\m{\mu}

\def\id{\operatorname{id}}
\def\FF{{\Bbb F}}

\topmatter

\title{Cluster algebras and Poisson geometry}
\endtitle
\author Michael Gekhtman$^*$, Michael Shapiro$^\ddag$,
 and Alek Vainshtein$^\dag$ \endauthor
\affil $^*$ Department of Mathematics, University of Notre Dame,
Notre Dame, IN 46556,
{\tt Michael.Gekhtman.1\@nd.edu}\\
$^\ddag$ Matematiska Institutionen, KTH, Stockholm,
{\tt mshapiro\@math.kth.se}\\
and Department of Mathematics, Michigan State University,
East Lansing, MI 48824-1027, {\tt mshapiro\@math.msu.edu}\\
$^\dag$ Departments of Mathematics and of Computer Science,
University of Haifa,  Israel 31905,
{\tt alek\@mathcs.haifa.ac.il}
\endaffil
\leftheadtext{Michael Gekhtman, Michael Shapiro, and Alek Vainshtein}

\abstract
We introduce a Poisson variety compatible with a cluster algebra
structure and a compatible toric action on this variety.
We study Poisson and topological properties of the union of generic
orbits of this toric action.
In particular, we compute the number of connected components of the union
of generic toric orbits for cluster algebras over real numbers.
As a corollary we compute the number of connected components
of refined open Bruhat cells in Grassmanians $G(k,n)$ over $\R$.
\endabstract
\endtopmatter

\heading{0. Introduction}
\endheading

The aim of the present paper is to study Poisson structures
naturally related to cluster algebras recently introduced
by Fomin and Zelevinsky in \cite{FZ2}.

Roughly speaking, a cluster algebra is defined by
an $n$-regular tree whose vertices correspond to $n$-tuples of cluster
variables  and edges describe birational transformations
between two $n$-tuples of variables; the cluster algebra itself is generated
by the union of all cluster variables.
Model examples of cluster algebras
are coordinate rings of double Bruhat cells
(see~\cite{FZ2}).
Edge transformation rules  imitate simplest (3-term) Pl\"ucker relations.
Given the set of transformation rules for all the edges incident to one vertex
of the tree, one can restore all the other transformation rules. The evolution
of transformation rules provides the so-called
Laurent phenomenon (\cite{FZ3}), which  means the following:
fix one cluster with cluster variables
$x_1,\dots,x_n$ and express any other cluster variable in terms of
$x_1,\dots,x_n$; then the expression is a Laurent
polynomial in $x_1,\dots,x_n$.

The first goal of this paper is to give less formal explanation for the
evolution of edge
transformation rules from the Poisson point of view.
Namely, we introduce a Poisson structure
compatible with the cluster algebra structure.
Compatibility simply means that the Poisson structure is homogeneously
quadratic in any set of cluster variables.
Then edge transformations describe simple transvections with respect to the
Poisson structure. In particular, the transition of transformation rules can
be simply explained as the transformation of coefficients of the Poisson
structure under a transvection.

The second goal of our paper is to extend our previous calculations of the
number of connected components in double Bruhat cells (see~\cite{SSV1, SSV2,
SSVZ, GSV})  to a more general setting of geometric cluster
algebras and compatible Poisson structures.
Namely, given a cluster algebra $\AA$ over $\R$
we compute the number of connected
components in the union
of generic symplectic leaves of any compatible Poisson structure in
a certain ``large'' nonsingular subset of $\Spec(\AA)$.

Finally, we apply the general formula to a special case of Grassmannian
coordinate ring.

The structure of the paper is as follows.

In the first chapter we
recall a notion of (geometric) cluster algebra $\AA$, introduce a notion
of a Poisson bracket compatible with $\AA$ and describe
all Poisson brackets compatible with $\AA$.
Moreover, we also prove the
following (partial) inverse result. Assume that a
homogeneous quadratic Poisson bracket on a rational $n$-dimensional manifold
is given (recall that the field of meromorphic
functions on such a manifold is a transcendental extension of the ground
field, i.e. the field of rational functions in $n$ variables).
We are looking for birational involutions preserving quadratic
homogeneity and satisfying some locality and universality properties.
Then there exists a cluster algebra compatible with the Poisson structure,
such that
these birational transformations are exactly edge transitions for this
cluster algebra.

In the second chapter we introduce an $\F^*$-action compatible with the
cluster
algebra $\AA$ (here $\F$ is a field of characteristic $0$).
Compatibility of the
$\F^*$-action means that all edge transformations
of $\AA$ are preserved under this action.
The union $\X^0$ of generic orbits with respect to this $\F^*$-action
is ``almost'' the union of generic symplectic leaves of a compatible Poisson
structure in $\Spec(\AA)$. We compute the number of connected components of
$\X^0$ for a cluster algebra over $\R$.

Finally we apply the previous result to the case of refined open Bruhat
cells in the Grassmannian $G(k,n)$.
Namely, the famous Sklyanin Poisson-Lie bracket on $SL_n(\R)$ induces a
Poisson bracket on the open Bruhat cell in  $G(k,n)$. This Poisson
bracket is compatible with the structure of a special cluster algebra,
one of whose clusters consists only of Pl\"ucker coordinates. The corresponding
$\R^*$-action determines the union of generic orbits, which is
simply described
as a subset of the Grassmannian defined by inequalities
$X_i\ne 0$, $i\in [1,n]$,
where $X_i$ is  the (cyclically solid) minor containing
the $i$th, $(i+1)$th, $\dots$,
$(i+k)$th ($\mod n$) columns.  We call this subset a refined open Bruhat
cell in the Grassmannian $G(k,n)$;
indeed, this subset is an intersection of $n$ open Bruhat cells in general
position.
In the last part we compute the number of connected components of a refined
open Bruhat cell in
$G(k,n)$ over $\R$.
This number is equal to $3\cdot 2^{n-1}$ if $k\gs 3$ and $n\gs 7$.

The authors would like to thank S.~Evans, S.~Fomin, P.~Foth,
B.~Shapiro, A.~Zelevinsky for many valuable discussions.

The first and the third authors are grateful to the Gustafsson foundation
for the financial support of their visits to KTH in the Fall 2000
and in the Spring 2001.
The second and the third authors express their gratitude to the Max Planck
Institute f\"ur Mathematik in Bonn,
where the final version of this paper was prepared.

\heading {1. Cluster algebras of geometric type  and Poisson brackets}
\endheading

\subheading {1.1. Cluster algebras of rational functions on a rational
$n$-dimensional manifold}
Let $A$ be an arbitrary matrix, $I=\{i_1,\dots,i_m\}$, $J=\{j_1,\dots,
j_n\}$ be two ordered multi-indices. We denote by $A(I;J)$ the $m\times n$
submatrix of $A$ whose entries lie in the rows $i_1,\dots,i_m$ and
columns $ j_1,\dots,j_n$. Instead of $A([1,m];[1,n])$ we write just $A[m;n]$
(here and in what follows we use the notation $[i,j]$ for a contiguous index
set $\{i,\ldots,j\}$).
Given a diagonal matrix
$D$ with positive integer diagonal entries $d_1,\dots,d_m$,
let $\ZZ_{mn}^D$ be the set of all $m\times n$ integer matrices $Z$
such that $m\ls n$ and $Z[m;m]$ is $D$-skew-symmetrizable (that is,
$DZ[m;m]$ is skew-symmetric); clearly, $\ZZ_{mn}^D=\ZZ_{mn}^{\lambda D}$
for any positive integer $\lambda$.
According to \cite{FZ2}, any $Z=(z_{ij})\in\ZZ_{mn}^D$
defines a cluster algebra of geometric type in the following way.
Let us fix a set of $m$ cluster variables $f_1,\dots,f_m$, and a set of
$n-m$ tropic variables $f_{m+1},\dots,f_n$. For each $i\in [1,m]$
we introduce a transformation $T_i$ of cluster variables by
$$
T_i(f_j)=\bar f_j=\left\{\alignedat2
&\frac1{f_i}\left(\prod_{z_{ik}>0} f_k^{z_{ik}}+\prod_{z_{ik}<0}
  f_k^{-z_{ik}}\right)\quad& &\text{for $j= i$}\\
&f_j \qquad& &\text{for $j\ne i$},
\endalignedat\right.\tag1.1
$$
and the corresponding matrix transformation $\bar Z=T_i(Z)$, called mutation,
by
$$
\bar z_{kl}=\left\{\alignedat2
&-z_{kl}\qquad& &\text{for $(k-i)(l-i)=0$}\\
&z_{kl}+\frac{|z_{ki}|z_{il}+z_{ki}|z_{il}|}2\quad& &\text{for $(k-i)(l-i)\ne
  0$}.
\endalignedat\right. \tag1.2
$$
Observe that the tropic variables are not affected by $T_i$, and that $\bar Z$
belongs to $\ZZ_{mn}^D$. Thus, one can apply transformations $T_i$ to the new
set of cluster variables (using the new matrix), etc. The cluster algebra
(of geometric type) is the subalgebra of the field of rational functions in
cluster variables $f_1,\dots,f_m$ generated by the union of all clusters;
its ground ring is the ring of integer polynomials over tropic variables.
We denote this algebra by $\AA(Z)$.

One can represent $\cC$ with the help of an $m$-regular tree $\Bbb T_m$
whose edges are labeled by the numbers $1,\dots,m$ so that the $m$ edges
incident to each vertex receive different labels. To each vertex $v$ of
$\Bbb T_m$ we assign a set of $m$ cluster variables $f_{v,1},\dots,f_{v,m}$
and a set of $n-m$ tropic variables $f_{m+1},\dots,f_n$. For an edge
$(v,\bar v)$ of $\Bbb T_m$ that is labeled by $i\in [1,m]$, the
variables $f=f_v$ and $\bar f=f_{\bar v}$ are related by the transformation
$T_i$ given by (1.1). The first monomial in the right hand side of (1.1)
is sometimes denoted by $M^i=M_v^i$. Transformations (1.2) then guarantee that
the second monomial is $\bar M^i=M_{\bar v}^i$.

Let us say that  cluster and tropic variables together form an extended
cluster. Assume that the entries of the initial extended cluster are
coordinate functions on a rational $n$-dimensional manifold $\Cal M^n$.
We thus get a realization of a cluster algebra
of geometric type as a cluster algebra of rational functions on $\Cal M^n$.
It is easy to see that in this situation
entries of any extended cluster are functionally independent.

\demo{Remark}
If all entries of $Z$ belong to the set $\{0, \pm 1\}$, it is sometimes
convenient to represent $Z$ by
a directed graph $\gr$ with vertices corresponding to the variables (both
cluster and tropic)
and with edges $i\to j$ for every pair of vertices $i,j$, such that
$Z_{ij}=1$ (in particular, there are no edges between vertices corresponding
to tropic variables).
If we assume, in addition, that the resulting graph has no nonoriented
3-cycles,
then the graph that corresponds to $\bar Z$ differs from the one that corresponds to $Z$
as follows. All edges through $i$ change directions. Furthermore, for every
two vertices $j,k$ such that edges $j\to i$ and $i \to k$ belong to the
graph of $Z$,  the graph of  $\bar Z$ contains an edge $j \to k$ if and only if
the graph of $Z$ does not contain an edge $k\to j$.
\enddemo

\subheading{1.2. $\tau$-coordinates}
For our further purposes it is convenient to consider, along with $f$, another
$n$-tuple of rational functions. In what follows this $n$-tuple is denoted
$\tau=(\tau_1,\dots,\tau_n)$, and is related to the initial
$n$-tuple $f$ as follows. Let $\widehat{D}$ be an $n\times n$ diagonal matrix $\widehat{D}=
\diag(d_1,\dots, d_m,1,\dots,1)$.
Denote by $\ZZ_{nn}^{\widehat{D}}(Z)$ the set of all
$n\times n$ integer matrices $Z'\in\ZZ_{nn}^{\widehat{D}}$ such that $Z'[m;n]=Z$.
Fix an arbitrary matrix $Z'\in\ZZ_{nn}^{\widehat{D}}(Z)$ and put
$$
\tau_j=f_j^{\varkappa_j}\prod_{k=1}^nf_k^{z'_{jk}}, \tag 1.3
$$
where $\varkappa_j$ is an integer, $\varkappa_j=0$ for $1\ls j\ls m$.
Given an extended cluster $f$, we say that the entries of the
corresponding $\tau$ form a $\tau$-cluster.

We say that the transformation $f\mapsto\tau$ is {\it nondegenerate\/}
if
$$
\det(Z'+K)\ne0,\tag1.4
$$
where $K=\diag(\varkappa_1,\dots,\varkappa_n)$. It is easy to see that
if the transformation $f\mapsto\tau$ is nondegenerate and the entries
of the extended cluster are functionally independent, then so are the
entries of the $\tau$-cluster.

\proclaim{Lemma 1.1} A nondegenerate transformation $f\mapsto\tau$ exists
if and only if $\rank Z=m$.
\endproclaim

\demo{Proof} The only if part is trivial. To prove the if part, assume that
$\rank Z=m$ and $\rank Z[m;m]=k\ls m$. Then there exists
a nonzero $m\times m$ minor of $Z$ contained in the columns $j_1<
j_2<\dots<j_m$ so that $j_k\ls m$, $j_{k+1}>m$ (here $j_0=0$, $j_{m+1}=n+1$).
 Without loss of generality assume that $j_i=i$
for $i\in[1,k]$. Define
$$
\varkappa_j=\left\{\alignedat2
&0\qquad& &\text{for $j=1,\dots,m$},\\
&1\qquad& &\text{for $j=j_{k+1},\dots,j_m$},\\
&\varkappa \qquad& &\text{otherwise}.
\endalignedat\right.
$$
Let us prove that there exists an integer $\varkappa$ such that
$\det(Z'+K)\ne0$. Indeed, the leading coefficient of this determinant
(regarded as a polynomial in $\varkappa$) is equal to the
$(2m-k)\times(2m-k)$ minor contained in the rows and columns $1,2,\dots,m,
j_{k+1},\dots,j_m$.
It is easy to see that using the same elementary row and
column operations one can reduce the corresponding submatrix to the form
$$
M=\pmatrix Z_1 & 0   & 0\\
               0   & Z_2 & Z_3\\
               0   & Z_4 & Z_5\endpmatrix,
$$
where $Z_1$ is just $Z[k;k]$,
$Z_2$ is an $(m-k)\times(m-k)$ matrix depending
only on the entries of $Z[m;m]$, $Z_3$, $Z_4$, $Z_5$ are $(m-k)\times
(m-k)$ matrices. Moreover, $Z_2=0$, since otherwise $\rank Z[m;m]$
would exceed $k$, and hence $\det M=\det Z_1\det Z_3\det Z_4$.
On the other hand, condition $\rank Z=m$ implies
$\det Z_1\det Z_3\ne0$, while the skew-symmetrizability of $Z'$ implies
$\det Z_4\ne0$. Therefore, the leading coefficient of $\det(Z'+K)$ is
distinct from zero, and we are done.
\qed
\enddemo

Let us find explicit expressions for the transformations $T_i$ in the
new coordinates $(\tau, Z')$. First, we extend the rules (1.2) to all
the entries of $Z'$. Observe that for any $Z'$ as above, one has
$\bar Z'[m;n]=\bar Z$. Moreover,
  the coordinate change $\bar f\mapsto\bar\tau$ remains nondegenerate,
  due to Lemma~1.1 and the following proposition.

\proclaim{Lemma 1.2} If $\rank Z=m$, then $\rank\bar Z=m$.
\endproclaim

\demo{Proof} Indeed, consider the following sequence of row and column
operations with the matrix $Z'$. For any $l$ such that $z'_{il}<0$,
subtract the $i$th column multiplied by $z'_{il}$ from the $l$th
column. For any $k$ such that $z'_{ki}>0$, add the $i$th row
multiplied by $z'_{ki}$ to the $k$th row. Finally, multiply the $i$th
row and column by $-1$. It is easy to see that the result of these
operations is exactly $\bar Z'$, and the lemma follows.
\qed
\enddemo

Finally, coordinates $\tau$ are transformed as follows.

\proclaim{Lemma 1.3} Let $i\in [1,m]$ and let
$\bar\tau_j=T_i(\tau_j)$ for $j\in [1,n]$. Then
$\bar\tau_i=1/\tau_i$ and $\bar\tau_j=\tau_j\psi_{ji}(\tau_i)$, where
$$
\psi_{ji}(\xi)=\left\{\alignedat2
\bigg(&\frac{1}{\xi}+1\bigg)^{-z'_{ji}}\qquad& &\text{for $z'_{ji}>0$},\\
(&{\xi}+1)^{-z'_{ji}}\qquad& &\text{for $z'_{ji}<0$},\\
&1 \qquad& &\text{for $z'_{ji}=0$ and  $j\ne i$}.
\endalignedat\right.
$$
\endproclaim

\demo{Proof} Let us start from the case $j=i$. Since $i\in [1,m]$ and
$z'_{ii}=0$, we can write
$$
\bar \tau_i=\prod_{k\ne i}\bar f_k^{\bar z'_{ik}}=\prod_{k\ne i}
f_k^{-z'_{ik}}=\frac{1}{\tau_i},
$$
as required.

Now, let $j\ne i$. Then
$$\multline
\bar\tau_j=\bar f_j^{\varkappa_j}
\bar f_i^{\bar z'_{ji}}\prod_{k\ne i}\bar f_k^{\bar z'_{jk}}\\=
 f_j^{\varkappa_j}
\left(\prod_{z'_{ik}>0}f_k^{z'_{ik}}+\prod_{z'_{ik}<0}
f_k^{-z'_{ik}}\right)
^{-z'_{ji}}f_i^{z'_{ji}}\prod_{k\ne i} f_k^{z'_{jk}}
\prod_{k\ne i}f_k^{(|z'_{ji}|z'_{ik}+z'_{ji}|z'_{ik}|)/2}\\=
\tau_j\left(\prod_{z'_{ik}>0}f_k^{z'_{ik}}+
\prod_{z'_{ik}<0}f_k^{-z'_{ik}}\right)
^{-z'_{ji}}\prod_{k\ne i}f_k^{(|z'_{ji}|z'_{ik}+z'_{ji}|z'_{ik}|)/2}.
\endmultline
$$

If $z'_{ji}=0$, then evidently $\bar\tau_j=\tau_j$.

Let $z'_{ji}>0$, then
$$\multline
\bar\tau_j=
\tau_j\left(\prod_{z'_{ik}>0}f_k^{z'_{ik}}+
\prod_{z'_{ik}<0}f_k^{-z'_{ik}}\right)
^{-z'_{ji}}\prod_{z'_{ik}>0}f_k^{z'_{ji}z'_{ik}}\\=\tau_j
\left(\prod_{z'_{ik}\ne 0}f_k^{-z'_{ik}}+1\right)^{-z'_{ji}}=
\tau_j(1/\tau_i+1)^{-z'_{ji}},
\endmultline
$$
as required.

Let  $z'_{ji}<0$, then
$$\multline
\bar\tau_j=
\tau_j\left(\prod_{z'_{ik}>0}f_k^{z'_{ik}}+
\prod_{z'_{ik}<0}f_k^{-z'_{ik}}\right)
^{-z'_{ji}}\prod_{z'_{ik}<0}f_k^{-z'_{ji}z'_{ik}}\\=\tau_j
\left(\prod_{z'_{ik}\ne 0}f_k^{z'_{ik}}+1\right)^{-z'_{ji}}
\prod_{z'_{ik}\ne 0}f_k^{-z'_{ij}z'_{ik}}=
\tau_j(\tau_i+1)^{-z'_{ji}},
\endmultline
$$
as required.
\qed
\enddemo

\subheading{1.3. Poisson brackets}
Let $\omega$ be a Poisson bracket on an $n$-dimensional manifold.
We say that functions $g_1,\dots,g_n$ are {\it log-canonical}
with respect to $\omega$ if $\omega(g_i,g_j)=\omega_{ij} g_ig_j$, where
$\omega_{ij}$ are integer constants.
The matrix $\Omega=(\omega_{ij})$ is called
the {\it coefficient matrix\/} of $\omega$ (in the basis $g$);
evidently, $\Omega\in so_n(\Z)$.

Fix some $Z\in\ZZ_{mn}^D$ and consider the following question:
are there any Poisson structures on a rational $n$-dimensional
manifold such that all clusters in the
cluster algebra $\AA(Z)$ are log-canonical with respect to them?

We say that a skew-symmetrizable matrix $A$ is {\it reducible\/} if there
exists a permutation matrix $P$ such that $PAP^T$ is a block-diagonal
matrix, and {\it irreducible\/} otherwise;
$r(A)$ is defined as the maximal number of diagonal blocks in $PAP^T$.
The partition into blocks defines an obvious equivalence relation
$\sim$ on the rows (or columns) of $A$.

\proclaim{Theorem 1.4}
Assume that $Z\in \ZZ_{mn}^D$ and $\rank Z=m$.
Then the Poisson brackets on a rational $n$-dimensional manifold
for which all extended clusters in $\AA(Z)$ are log-canonical
form a vector space of
dimension $r+\binom{n-m}2$, where $r=r(Z[m;m])$.
Moreover, the coefficient matrices of these Poisson brackets
in the basis $\tau$ are characterized by the
equation $\Omega^\tau[m;n]=\Lambda Z\widehat{D}^{-1}$,
where $\Lambda=\diag(\lambda_1,\dots,\lambda_m)$ with
$\lambda_i=\lambda_j$ whenever $i\sim j$.
In particular, if $Z[m;m]$ is irreducible, then
$\Omega^\tau[m;n]=\lambda Z\widehat{D}^{-1}$.
\endproclaim

\demo{Proof} Let us note first that
$\tau$-coordinates are expressed in a monomial way in terms of initial
coordinates $f$, and that this transformation is invertible.
Therefore, all extended clusters in $\AA(Z)$ are log-canonical w.r.t. some
bracket $\omega$ if and only if so are all the
 corresponding $\tau$-clusters. Denote by $\Omega^f$ and by $\Omega^\tau$
the matrices of $\omega$ in the bases $f$ and $\tau$, respectively. It
is easy to see that $\Omega^\tau=(Z'+K)\Omega^f(Z'+K)^T$. Evidently,
transformation $T_i$ preserves the log-canonicity if and only if for any
$j\ne i$,
$\omega(\bar f_i,\bar f_j)=\bar\omega_{ij}\bar f_i\bar f_j$ provided
$\omega(f_i,f_j)=\omega_{ij}f_if_j$. Using (1.1) we get
$$
\multline
\omega(\bar f_i,\bar f_j)=\omega\left(\frac1{f_i}\left(\prod_{
z_{ik}>0}f_k^{z_{ik}}+\prod_{z_{ik}<0}f_k^{-z_{ik}}\right),f_j\right)\\
=\frac{f_j}{f_i}\prod_{z_{ik}>0}f_k^{z_{ik}}
\left(\sum_{z_{ik}>0}z_{ik}\omega_{kj}-\omega_{ij}\right)+
\frac{f_j}{f_i}\prod_{z_{ik}<0}f_k^{z_{ik}}
\left(-\sum_{z_{ik}<0}z_{ik}\omega_{kj}-\omega_{ij}\right),
\endmultline
$$
and hence the above conditions are satisfied if and only if
$\sum_{z_{ik}>0}z_{ik}\omega_{kj}-\omega_{ij}=
-\sum_{z_{ik}<0}z_{ik}\omega_{kj}-\omega_{ij}$ for $j\ne i$.
This means that
$$
(Z'+K)\Omega^f[m;n]=Z\Omega^f=\pmatrix \Delta &0\endpmatrix,
\tag 1.7
$$
where $\Delta$ is a diagonal matrix. Consequently, we get
$\Omega^\tau[m;n]=\Delta (Z'[n;m])^T$, and hence
$\Delta Z^T[m;m]=\Omega^\tau[m;m]$ is skew-symmetric.
Therefore, $\Delta=-\Lambda D^{-1}$
where $\Lambda=\diag(\lambda_1,\dots,\lambda_m)$ with
$\lambda_i=\lambda_j$ whenever $i\sim j$.
It remains to notice that $Z'=-\D^{-1}{Z'}^T\D$,
and therefore $Z'[n;m]=-\D^{-1}Z^TD$, and the equation
$\Omega^\tau[m;n]=\Lambda Z\D^{-1}$ follows. The entries $\omega^f_{ij}$
for $m+1\ls i<j\ls n$ are free parameters.
\qed\enddemo

\subheading{1.4. Recovering cluster algebra transformations}
In this section we recover transformations (1.1), (1.2) (in their
equivalent form presented in Lemma~1.3) as unique involutive transformations
of log-canonical bases satisfying certain additional restrictions.

By the definition, {\it local data\/} $F$ is
a family of rational functions in one variable $\psi_w$,
$w=0,\pm1,\pm2,\dots$, and an additional function in one
variable $\fy$.
For any Poisson bracket $\omega$ and any log-canonical (with respect
to $\omega$) basis $t=(t_1,\dots,t_n)$, the local data $F$ gives rise
to $n$ transformations $F_i^\omega$ defined as follows:

(i) $F_i^\omega(t_i)=\bar t_i=\fy(t_i)$;

(ii) let $\Omega=(\omega_{ij})$  be the coefficient matrix of $\omega$ in
the basis $t$, then $F_i^\omega(t_j)=\bar t_j=t_j\psi_{\omega_{ij}}(t_i)$
for $j\ne i$.

We say that  local data $F$
is {\it canonical\/} if
for any Poisson bracket $\omega$, any log-canonical (with respect
to $\omega$) basis $t$, and any index $i$, the set $F_i^\omega(t)$ is
a log-canonical basis of $\omega$ as well. Local data is called
{\it involutive\/} if any $F_i^\omega$ is an involution, and
is called {\it normalized\/} if $\psi_w(0)=\pm 1$
for any integer $w\gs0$.

We say that a polynomial $P$ of degree $p$ is {\it $a$-reciprocal\/}
if $P(0)\ne0$ and there
exists a constant $c$ such that $\xi^pP(a/\xi)=cP(\xi)$ for any
$\xi\ne0$.

The following result gives a complete description of normalized
involutive canonical local data.

\proclaim{Lemma 1.5} Any normalized involutive canonical local
data has one of the following forms:

{\rm (i)} $\fy(\xi)=\xi$ and $\psi_w(\xi)=\pm1$ for any integer $w$
{\rm(}trivial local data{\rm);}

{\rm (ii)} $\fy(\xi)=-\xi$ and
$\psi_w(\xi)=\pm\left(\frac{P(\xi)}{P(-\xi)}\right)^w$, where $P$ is a
polynomial without symmetric roots;

{\rm (iii)} $\fy(\xi)=\frac a{\xi}$,
$\psi_w(\xi)=a_w\xi^{c_w}\psi_1^w(\xi)$, and
$\psi_1(\xi)=\frac{P(\xi)}{Q(\xi)}$, where $P$ and $Q$ are coprime
$a$-reciprocal polynomials of degrees $p$ and $q$,
and the constants $a_w$, $c_w$, $p$, $q$ satisfy relations
$a_{-1}^2=a^{-c_{-1}}$,  $p-q=c_{-1}$, and
$$
a_w=\left\{\alignedat2
&\pm1\qquad& &\text{for $w\gs0$},\\
&a_{-w}^{-1}a_{-1}^{-w}\quad& &\text{for $w<0$},
\endalignedat\right.\qquad
c_w=\left\{\alignedat2
&0\qquad& &\text{for $w\gs0$},\\
&-wc_{-1}\quad& &\text{for $w<0$}.
\endalignedat\right.
$$
\endproclaim

\demo{Proof} Let $F$ be canonical local data; fix an arbitrary $\omega$
and pick any nonzero entry $\omega_{ij}$. Applying $F_i^\omega$ gives
$$
\omega(\bar t_i,\bar t_j)=\omega(\fy(t_i),\psi_w(t_i)t_j)=
\frac{\fy'(t_i)t_i}{\fy(t_i)}w\bar t_i\bar t_j,
$$
where $w=\omega_{ij}$.
The canonicity of $F$ yields
$$
\frac{\fy'(\xi)\xi}{\fy(\xi)}=c,
$$
where $cw$ is an integer, therefore $\fy(\xi)=a\xi^c$. Since $F$ is
involutive, we get $a^{c+1}\xi^{c^2}=\xi$, hence either $c=1$ and
$a=\pm1$, or $c=-1$ and $a$ is an arbitrary nonzero constant. Let
$\bar\Omega=(\bar\omega_{ij})$ be the coefficient matrix of $\omega$
in the basis $\bar t$. It follows immediately from the above
calculations that $\bar \omega_{ij}=\omega_{ij}$ if $c=1$ and
$\bar\omega_{ij}= -\omega_{ij}$ if $c=-1$.

Assume that $c=1$ and $a=1$. Then the involutivity of $F$ (applied to
$t_j$) gives $\psi^2_w(\xi)\equiv1$ for any $w$, and hence
$\psi_w(\xi)=\pm1$, and we obtained the first case described in the
lemma.

To proceed further, consider two arbitrary indices $j,k\ne i$ and put
$u=\omega_{ij}$, $v=\omega_{ik}$. Then we get
$$
\omega(\bar t_j,\bar t_k)=\omega(t_j\psi_u(t_i),t_k\psi_v(t_i))=
\left(\omega_{jk}+v\frac{\psi_u'(t_i)t_i}{\psi_u(t_i)}-
u\frac{\psi_v'(t_i)t_i}{\psi_v(t_i)}\right)\bar t_j\bar t_k,
$$
so the canonicity of $F$ yields
$$
v\frac{\psi_u'(\xi)\xi}{\psi_u(\xi)}-
u\frac{\psi_v'(\xi)\xi}{\psi_v(\xi)}=c_{uv}
$$
for some integer constant $c_{uv}$. Integrating both sides we obtain
equation $\psi_u^v(\xi)=a_{uv}\xi^{c_{uv}}\psi_v^u(\xi)$, which is
valid for any $u$ and $v$. In particular, taking $v=1$ we get
$$
\psi_u(\xi)=a_u\xi^{c_u}\psi_1^u(\xi),\tag 1.8
$$
where $a_u=a_{u1}$ and $c_u=c_{u1}$.

Let us return to the case $\fy(\xi)=-\xi$. In this case the
involutivity of $F$ applied to $t_j$ gives
$\psi_w(\xi)\psi_w(-\xi)\equiv1$ for any integer $w$. Using (1.8) we
get
$$
(-1)^{c_w}a_w^2\xi^{2c_w}\bigl(\psi_1(\xi)\psi_1(-\xi)\bigr)^w\equiv1,
$$
which immediately yields $c_w=0$ and $a_w=\pm1$.
Let us represent the
rational function $\psi_1$ as the ratio of two coprime polynomials
$P$ and $Q$. Then the above involutivity condition gives
$P(\xi)P(-\xi)=Q(\xi)Q(-\xi)$, which means that $Q(\xi)=\pm P(-\xi)$.
Therefore,
$$
\psi_w(\xi)=\pm\left(\frac{P(\xi)}{P(-\xi)}\right)^w,
$$
and the coprimality condition translates to the nonexistence of
symmetric roots of $P$.

Finally, consider the case $\fy(\xi)=a/\xi$. The involutivity of $F$
applied to $t_j$ gives $\psi_w(\xi)\psi_{-w}(a/\xi)\equiv1$ for any
integer $w$. Using (1.8) we get
$$
\left(\frac{\psi_1(a/\xi)}{\psi_1(\xi)}\right)^w=
a^{c_{-w}}a_wa_{-w}\xi^{c_w-c_{-w}};\tag 1.9
$$
in particular, for $w=-1$ this can be rewritten as
$\psi_1(\xi)=a_{-1}\xi^{c_{-1}}\psi_1(a/\xi)$, and for $w=1$, as
$\psi_1(a/\xi)=a_{-1}a^{c_{-1}}\xi^{-c_{-1}}\psi_1(\xi)$ (since by definition,
$a_1=1$ and $c_1=0$). Comparing the two latter identities
one immediately gets  $a_{-1}^2=a^{-c_{-1}}$. Besides, plugging
$$
\frac{\psi_1(a/\xi)}{\psi_1(\xi)}=a_{-1}a^{c_{-1}}x^{-c_{-1}}
$$
into (1.9), one gets $c_w-c_{-w}=-wc_{-1}$ and $a_wa_{-w}=
a_{-1}^wa^{wc_{-1}-c_{-w}}=a_{-1}^wa^{-c_w}$.  Observe that the normalization
condition applied to (1.8) immediately
gives $c_w=0$ and $a_w=\pm1$ for $w\gs0$. Therefore, we get
$c_w=-wc_{-1}$ and $a_w=a_{-w}^{-1}a_{-1}^{-w}$ for $w<0$.

Finally,
let us represent the rational function $\psi_1$ as
$\psi_1(\xi)=\frac{P(\xi)}{Q(\xi)}$, where $P$ and $Q$ are
coprime polynomials; by the normalization condition, $P(0)/Q(0)=\pm1$.
Plugging this into (1.9) for $w=1$ one gets
$$
a_{-1}\xi^{c_{-1}}\frac{P(a/\xi)}{Q(a/\xi)}=\frac{P(\xi)}{Q(\xi)}.
$$
Since  $P$ and $Q$ are coprime, the above identity can only hold when
they both are $a$-reciprocal.
Equating the degrees on both sides of the identity gives $p-q=c_{-1}$.
\qed
\enddemo

We say that local data $F$ is {\it finite\/} if it
has the following finiteness property: let $n=2$, and let $\omega$
possess a log-canonical basis $t=(t_1,t_2)$ such that the corresponding
coefficient matrix  has the simplest form
$\pmatrix 0 & 1\\-1 &0\endpmatrix$; then the group generated by
$F_1^\omega$  and $F_2^\omega$ has a finite order.

\proclaim{Theorem 1.6} Any nontrivial finite involutive canonical local data
has one of the following forms:

{\rm (i)} $\fy(\xi)=a/\xi$ with $a\ne0$, $\psi_w(\xi)=(\pm1)^wa_w$, where
$a_w=\pm1$ and $a_w=a_{-w}$;

{\rm (ii)} $\fy(\xi)=b^2/\xi$ with $b\ne0$,
$$
\psi_w(\xi)=\left\{\alignedat2
&(\pm1)^wa_w\left(\frac{\xi+b}{b}\right)^w\qquad& &\text{for $w\gs0$},\\
&(\pm1)^wa_{-w}\left(\frac{\xi+b}{\xi}\right)^w\qquad& &\text{for $w<0$},
\endalignedat\right.
$$
where $a_w=\pm1$ for $w\gs0$.
\endproclaim

\demo{Proof}  Assume first that local data is of type (ii) described in
Lemma~1.5. Let us prove that the transformation of the plane $(x,y)$
corresponding  to the composition $T=F_2^\omega\circ F_1^\omega$ has an
infinite order. Indeed, let as start from the pair $x=R_1(\xi)/S_1(\xi)$,
$y=R_2(\xi)/S_2(\xi)$, where $R_i$ and $S_i$ are coprime, $i=1,2$, and
$$
s_2< r_2< s_1< r_1 \tag 1.10
$$
(here and in what follows a small letter denotes
the degree of the polynomial denoted by the corresponding capital letter,
e.g., $s_2=\deg S_2$). Let $T(x,y)=(\widehat{x},\widehat{y})$, then
$$
\widehat{y}=\mp\frac{R_2\widetilde P(R_1,S_1)}{S_2\widetilde P(-R_1,S_1))},
$$
where $\widetilde P(\xi,\zeta)=\zeta^pP(\xi/\zeta)$. Observe that if $R$
and $S$ are coprime and $P$ does not have symmetric roots then
$\widetilde P(R,S)$ and $\widetilde P(-R,S)$ are coprime as well. Therefore,
if  $\mp\widehat{R}_2/\widehat{S}_2$ is the representation of $\widehat{y}$
as the ratio of coprime polynomials, then
$\widehat{r}_2\gs pr_1-s_2$. Besides, it is easy to see that
$\widehat{r}_2-\widehat{s}_2=r_2-s_2>0$. Next,
$$
\widehat{x}=\mp\frac{R_1\widetilde P(\mp \widehat{R}_2,\widehat{S}_2)}
{S_1\widetilde P(\pm \widehat{R}_2,\widehat{S}_2)},
$$
so if  $\mp\widehat{R}_1/\widehat{S}_1$ is the representation of
$\widehat{x}$ as the ratio of coprime polynomials, then
$\widehat{s}_1\gs p\widehat{r}_2-r_1$. Besides, $\widehat{r}_1-\widehat{s}_1=
r_1-s_1>0$. Finally, $\widehat{s}_1-\widehat{r}_2\gs p\widehat{r}_2-r_1-
\widehat{r}_2\gs(p-1)(pr_1-s_2)-r_1>(p-1)(pr_1-r_1)-r_1=p(p-2)r_1$.
Therefore, relation (1.10) is preserved under the action of $T$, provided
$p\gs2$. It remains to notice that $\widehat{r}_2-r_2\gs pr_1-s_2-r_2$, and
hence $p\gs2$ implies $\widehat{r}_2>r_2$. Therefore, $r_2$ grows
monotonically with the iterations of $T$, and hence $T$ has an infinite order.

It remains to consider the case $p=1$. Without loss of generality we may
assume that $P(\xi)=\xi+b$, where $b\ne0$. The choice of the sign in
the expressions for $F_1^\omega$ and $F_2^\omega$ leads to the following four
possibilities:
$$\align
T(x,y) &=\left(x\frac{-b^2+bx+by+xy}{b^2-bx+by+xy},-y\frac{b+x}{b-x}\right),\\
T(x,y) &=\left(x\frac{b^2-bx-by-xy}{b^2-bx+by+xy},-y\frac{b+x}{b-x}\right),\\
T(x,y) &=\left(x\frac{b^2-bx+by+xy}{-b^2+bx+by+xy},y\frac{b+x}{b-x}\right),\\
T(x,y) &=\left(x\frac{b^2-bx+by+xy}{b^2-bx-by-xy},y\frac{b+x}{b-x}\right).
\endalign
$$
In the first case,
$$
T^{2k}(x,y)=\left(x-\frac{4k}{b^2}x^2y,
y+\frac{4k}{b^2}xy^2\right)+o((x+y)^3),
$$
in the second case,
$$
T^{2k}(x,y)=\left(x, y-\frac{4k}{b}xy\right)+o((x+y)^2),
$$
in the third case,
$$
T^{2k}(x,y)=\left(x+\frac{4k}{b}xy,y\right)+o((x+y)^2),
$$
and in the fourth case,
$$
T^k(x,y)=\left(x+\frac{2k}{b}xy, y+\frac{2k}{b}xy\right)+o((x+y)^2).
$$
Therefore, in all these cases $T$ has an infinite order.

Assume now that local data is of type (iii) described in Lemma~1.5.
Instead of looking at $T=F_2^\omega\circ F_1^\omega$ we are going to study
$\widehat{T}=\sigma\circ F_1^\omega$, where $\sigma(x,y)=(y,x)$. It is easy
to check that $T=\widehat{T}^2$, so $T$ has a finite order if and only if
$\widehat{T}$ has a finite order.

Consider a pair $x=R_1(\xi)/S_1(\xi)$,
$y=R_2(\xi)/S_2(\xi)$, where $R_i$ and $S_i$ are coprime, $i=1,2$,
and denote $\widehat{T}(x,y)=(\widehat{x},\widehat{y})$. Then
$$
(\widehat{x},\widehat{y})=\left(\frac{R_2\widetilde P(R_1,S_1)}
{S_2\widetilde Q(R_1,S_1)S_1^{p-q}},\frac{aS_1}{R_1}\right),
$$
where $\widetilde P$ and $\widetilde Q$ are defined as before.
Let $\widehat{R}_1/\widehat{S}_1$ and $\widehat{R}_2/\widehat{S}_2$ be
the representations of $\widehat{x}$ and $\widehat{y}$
as the ratios of coprime polynomials, then
$$\alignat2
&\widehat{r}_1\gs p\max\{r_1,s_1\}-s_2,\qquad&
& \widehat{r}_2=s_1,\\
&\widehat{s}_1\gs q\max\{r_1,s_1\}+s_1\max\{0,p-q\}-r_2,\qquad&
&\widehat{s}_2=r_1.
\endalignat
$$

Assume first that $p\gs2$. Then we start from a pair
$(x,y)$ satisfying an additional condition $r_1>s_2$. Observe that
$\widehat{r}_1-\widehat{s}_2\gs  p r_1-s_2-r_1\gs
(r_1-s_2) +(p-2)r_1>0$, and hence the
above additional condition is preserved under the action of $\widehat{T}$.
Moreover, $\widehat{s}_2=r_1>s_2$, which means that
$s_2$ grows monotonically with the iterations of $\widehat{T}$, and hence
$\widehat{T}$ has an infinite order.

Assume now that $p<2$ and $q\gs2$. Then we start from a pair
$(x,y)$ satisfying an additional condition $s_1>r_2$. Observe that
$\widehat{s}_1-\widehat{r}_2\gs q s_1-r_2-s_1\gs
(s_1-r_2) +(q-2)s_1>0$, and hence the
above additional condition is preserved under the action of $\widehat{T}$.
Moreover, $\widehat{r}_2=s_1>r_2$, which means that
$r_2$ grows monotonically with the iterations of $\widehat{T}$, and hence
$\widehat{T}$ has an infinite order.

It remains to consider the case $\max\{p,q\}<2$, which amounts to
the following four possibilities: $p=1$, $q=1$; $p=1$, $q=0$;
$p=0$, $q=1$; $p=0$, $q=0$.

Consider the first possibility, when $p=1$, $q=1$.
Recall that by Lemma~1.5,
both $P$ and $Q$ are $a$-reciprocal. It is easy to check that a linear
$a$-reciprocal polynomial can be represented as $c(\xi+b)$ with $b^2=a$
and $c\ne0$.
Taking into account the normalization condition, we can write
$\widehat{T}$ as follows: $\widehat{T}\:(x,y)\mapsto
(\pm y(x+b)/(x-b),b^2/x)$.
Consider the transformation $\widetilde{T}=\widehat{T}^8$.
In a neighborhood of the point $(0,0)$
$\widetilde{T}$ can be expanded as
$$
\widetilde{T}\:(x,y)\mapsto\left(x-\frac{8}{b^2} x^2y,
y+\frac{8}{b^2} xy^2\right)+o((x+y)^3).
$$
Therefore,
$\widetilde{T}^k\:(x,y)
\mapsto (x-8kx^2y/b^2, y+8kxy^2/b^2)+o((x+y)^3)$, and
hence $\widetilde{T}$ has an infinite order.

Consider the second possibility, when  $p=1$, $q=0$.
In this case the transformation $\widehat{T}$ can be written as
follows: $\widehat{T}\:(x,y)\mapsto (\pm y(x+b)/b,a/x)$ with $b^2=a$.
It is easy to check that in this case indeed $\widehat{T}^5=\id$.
Besides, by Lemma~1.5, $c_{-1}=p-q=1$ and $a_{-1}^2=b^{-2}$.
Therefore,
$$
\psi_w(\xi)=\left\{\alignedat2
&(\pm1)^wa_w\left(\frac{\xi+b}{b}\right)^w\qquad& &\text{for $w\gs0$},\\
&\frac{(\pm1)^w}{a_{-w}}\left(\frac{\xi+b}{\xi}\right)^w\qquad&
&\text{for $w<0$},
\endalignedat\right.
$$
where $a_w=\pm1$ for $w\gs0$.

Consider the third possibility, when  $p=0$, $q=1$.
In this case the transformation $\widehat{T}$ can be written as follows:
$\widehat{T}\:(x,y)\mapsto (cy/(x+b),a/x)$ with $b^2=a$ and $c=\pm b$.
Consider the transformation $\widetilde{T}=\widehat{T}^7$.
In a neighborhood of the  point $(\infty,0)$
$\widetilde{T}$ can be expanded in local coordinates $z=1/x$, $y$ as
$$\multline
\widetilde{T}\:(z,y)\mapsto\left(\frac {cz}b-\frac{c-b}{b^2} zy
-\frac{b^2-3bc+3c^2}{b^2}z^2y-\frac{b-c}{b^3} zy^2,\right.\\
\left.\frac{by}c+\frac{b(c-b)}czy-\frac{b^2(c-b)}cz^2y-
\frac{b-2c}czy^2\right)+o((z+y)^3).
\endmultline
$$
Therefore, $\widetilde{T}$ can have a finite order only if $c=b$, in which
case we get
$\widetilde{T}\:(z,y)\mapsto (z-z^2y, y+zy^2)+o((z+y)^3)$.
 In its turn, this gives $\widetilde{T}^k\:(z,y)
\mapsto (z-kz^2y, y+kzy^2)+o((z+y)^3)$, and
hence $\widetilde{T}$ has an infinite order.

Finally, consider the last possibility, when $p=q=0$.
In this case the transformation $\widehat{T}$ can be written as follows:
$\widehat{T}\:(x,y)\mapsto (\pm y,a/x)$.
It is easy to check that in this case indeed $\widehat{T}^4=\id$.
Besides, by Lemma~1.5, $c_{-1}=p-q=0$ and $a_{-1}^2=1$.
Therefore,
$$
\psi_w(\xi)=\left\{\alignedat2
&(\pm1)^wa_w\qquad& &\text{for $w\gs0$},\\
&\frac{(\pm1)^w}{a_{-w}}\qquad& &\text{for $w<0$},
\endalignedat\right.
$$
where $a_w=\pm1$ for $w\gs0$.
\qed
\enddemo

Observe that in the skew-symmetric case the transformations of type (ii)
described in Theorem~1.6 coincide up to a sign with those found in
Lemma~1.3. Indeed, in this case $z'_{ji}=-z_{ij}$, and hence Lemma~1.3
gives $\fy(\xi)=1/\xi$ and
$$
\psi_{ji}(\xi)=\psi_{z_{ij}}(\xi)=\left\{\alignedat2
(&\xi+1)^{z_{ij}}\qquad& &\text{for $z_{ij}\gs0$},\\
\bigg(&\frac{\xi+1}{\xi}\bigg)^{z_{ij}}\qquad& &\text{for $z_{ij}<0$}.
\endalignedat\right.
$$
It is easy to check that after the scaling $\tau\mapsto b\tau$ one gets
$\fy(\xi)=b^2/\xi$ and
$$
\psi_{ji}(\xi)=\psi_{z_{ij}}(\xi)=\left\{\alignedat2
&\left(\frac{\xi+b}{b}\right)^{z_{ij}}\qquad& &\text{for $z_{ij}\gs0$},\\
&\left(\frac{\xi+b}{\xi}\right)^{z_{ij}}\qquad& &\text{for $z_{ij}<0$}.
\endalignedat\right.
$$

\heading{2. The cluster manifold}
\endheading

In this section we construct an algebraic variety $\X$ (which we call
the {\it cluster manifold\/})
related to the cluster algebra $\AA(Z)$. Our approach is suggested by
considering coordinate rings of double Bruhat cells, which provide
main examples of cluster algebras.

\subheading{2.1. Poisson brackets on the cluster manifold}
It would be very natural to define the cluster manifold as
$\Spec(\AA(Z))$, since
$\Spec(\AA(Z))$ is the maximal manifold $M$ satisfying the following two
universal conditions:

a) all the cluster functions are regular functions on $M$;

b) for any pair $x_1$, $x_2$ of two distinct points on $M$ there exists a
cluster function $f\in \AA(Z)$ such that $f(x_1)\ne f(x_2)$.

However, it was shown by S.~Fomin that the Markov cluster algebra defined
in \cite{FZ2} is not finitely generated.
This observation shows that $\Spec(\cC)$ might be a rather complicated object.
Therefore we define $\X$ as a ``handy'' nonsingular
part of $\Spec(\cC)$.

We will describe $\X$ by means of charts and transition
functions. Assume that $\cC$ is given by an $m$-regular tree
$\Bbb T_m$ (see Section 1.1).
For each vertex $v$ of $\Bbb T_m$ we define the chart, that is, an open subset
$X=X_v\subset \X$ by
$$
X_v=\Spec(\F[f_{v,1},f_{v,1}^{-1},\dots,f_{v,m},
f_{v,m}^{-1},f_{m+1},\dots,f_{n}])
$$
(as before, $\F$ is a field of characteristic $0$).
An edge $(v,\bar v)$ of $\Bbb T_m$ labeled by a number $i\in [1,m]$
defines a transition function $X_v\to X_{\bar v}=\bar X$
by the equations $\bar f_{j}=f_{j}$ if $j\ne i$,
and the {\it three-term relation\/}
$\bar f_{i}f_{i}=M^i_v+M^i_{\bar v}$,
where $M^i_u$, $u\in \{v,\bar v\}$,
are some monomials in $f_1,\dots,f_n$.

Note that the tree $\Bbb T_m$ is connected, hence any pair of its vertices is
connected by a unique path. Therefore, the transition map between the charts
corresponding to two arbitrary vertices  can be computed as the composition
of the transitions along this path. Finally, put
$$
\X=\cup_{v\in {\Bbb T_m}} X_v.
$$

It follows immediately from the definition that $\X\subset \Spec(\cC)$.
However, $\X$ contains only such points $x\in \Spec(\cC)$
for which there exists a vertex $v$ of $\Bbb T_m$ whose cluster variables
form a coordinate system in some neighbourhood $\X(x)\subset \X$ of this point.

\demo{Example} Consider a cluster algebra $\AA_1$ over $\C$ given
by two clusters $\{f_{1},f_2,f_3\}$ and $\{\bar f_{1},f_2,f_3\}$
subject to one relation: $f_{1}\bar f_{1}=f_2^2+f_3^2$.

In this case $\Spec(\AA_1)=\Spec(\C[x,y,u,v]/\{xy-u^2-v^2=0\})$ is a
singular affine hypersurface $H\subset \C^4$ given by the equation
$xy=u^2+v^2$ and containing a singular point $x=y=u=v=0$. On the other hand,
$\X=H\setminus\{x=y=u^2+v^2=0\}$ is nonsingular.
\enddemo

In the general case, the following proposition stems immediately from
the above definitions and Theorem~1.4.

\proclaim{Lemma 2.1} The cluster manifold $\X$ is a nonsingular
rational manifold and
possesses a Poisson bracket such that for each vertex $v$ of $\Bbb T_m$
the corresponding extended cluster is log-canonical w.r.t. this bracket.
\endproclaim

Let $\omega$ be one of these Poisson brackets. Recall that a
{\it Casimir element\/} corresponding to $\omega$ is a function that is in
involution with all the
other functions on $\X$. All rational Casimir functions form a subfield
$\F_C$ in
the field of rational functions $\F(\X)$. The following proposition
provides a complete description of $\F_C$.

\proclaim{Lemma 2.2} $\F_C=\F(\m_1,\dots,\m_s)$, where $\m_j$ has a
monomial form
$$
\m_j=\prod_{i=m+1}^n f_i^{\alpha_{ji}}
$$
for some integer $\alpha_{ji}$, and $s=\corank\omega$.
\endproclaim

\demo{Proof} Fix a vertex $v$ of $\Bbb T_m$ and consider
an open subset
$$
X_v^0=X^0=\{x\in X_v : f_i(x)\ne 0, i\in [m+1, n]\}
$$
(from now on we omit in the notation
the dependence of $X$, $f$ and other objects on $v$).  Define the
$\tau$-coordinates as in Section 1.2 and note that each $\tau_i$
is distinct from~$0$ and from $\infty$ in $X^0$. Therefore,
$\widetilde\tau_i=\log\tau_{i}$ form a coordinate system in $X^{0}$, and
$\omega(\widetilde\tau_p,\widetilde\tau_q)=\omega^\tau_{pq}$, where
$\Omega^\tau=
(\omega^\tau_{pq})$ is the coefficient matrix of $\omega$ in the basis
$\tau$. The Casimir functions of $\omega$ that are linear in $\{\widetilde\tau_{i}\}$
 are given by the left nullspace $N_{\el}(\Omega^\tau)$ in the following way.
Since $\Omega^\tau$ is an integer matrix,
its left nullspace contains an integral lattice $L$. For
any vector $u=(u_1,\dots,u_n)\in L$, the sum $\sum_i
u_i\widetilde\tau_{i}$ is in involution with all the coordinates
$\widetilde\tau_{j}$. Hence the product
$\prod_{i=1}^n\tau_{i}^{u_i}$ belongs to $\F_C$; moreover $\F_C$ is
generated by the monomials $\m^{u}=\prod_{i=1}^n\tau_{i}^{u_i}$ for $s$
distinct vectors ${u\in L}$.

Let us calculate $\log\m^{u}=u\widetilde\tau^T$, where $\tilde\tau=
(\tilde\tau_1,\dots,\tilde\tau_n)$.
Recall that by (1.3), $\widetilde\tau^T=(Z'+K)\widetilde f^T$, where
$K$ is a diagonal matrix whose first $m$ diagonal entries are equal
to zero, and $\widetilde f=(\log f_1,\dots,\log f_n)$. So
$u\widetilde\tau^T=\alpha\widetilde f^T$, where $\alpha=u(Z'+K)$.
Further, consider the decompositions
$$
Z'=\pmatrix Z_1 & Z_2\\
            -Z_2^TD & Z_4\endpmatrix,\qquad
\Omega^\tau=\pmatrix \Omega_1 & \Omega_2\\
                     -\Omega_2^T & \Omega_4\endpmatrix,
$$
where $Z_1=Z'[m;m]$, $\Omega_1=\Omega^\tau[m;m]$. By Theorem~1.4,
$Z_1=\Lambda^{-1}\Omega_1D$ and $Z_2=\Lambda^{-1}\Omega_2$; moreover
it is easy to check that $\Lambda$ and $\Omega_1$ commute.
Let $u=\pmatrix u^1 &u^2\endpmatrix$ and $\alpha=\pmatrix
\alpha^1 &\alpha^2\endpmatrix$ be the
corresponding decompositions of $u$ and $\alpha$.
Since $u\Omega^\tau=0$, one has $u^1\Omega_1-u^2\Omega_2^T=0$.
Therefore $\alpha^1=u^1Z_1-u^2Z_2^TD=(u^1\Omega_1-u^2\Omega_2^T)\Lambda^{-1}D
=0$, and hence $\m^u=\prod_{i=m+1}^nf_i^{\alpha_i}$. Finally, $X^0$ is open and
dense in $X$, so the above defined rational Casimir functions can be extended to the
whole $X$, and hence to $\X$.
\qed
\enddemo

\subheading{2.2. Toric action on the cluster algebra}
Assume that an integer weight $w_v=(w_{v,1},\dots,w_{v,n})$ is given at any
vertex $v$ of the tree $\Bbb T_m$. Then we define
a {\it local toric action\/} on the cluster at $v$ as the
$\F^\star$-action given by the
formula $\{f_{v,1},\dots,f_{v,n}\}\mapsto \{f_{v,1}\cdot t^{w_{v,1}},\dots,
f_{v,n}\cdot t^{w_{v,n}}\}$. We say that local toric actions are compatible
if for any two extended clusters $C_v$ and $C_u$ the following diagram is
commutative:
$$
\CD
C_v @>>> C_u\\
@V t^{w_v}VV @V t^{w_u}VV\\
C_v @>>> C_u
\endCD
$$
In this case, local toric actions together define a
global toric action on $\cC$. This toric action is said to be an
{\it extension\/} of the above local actions.

A toric action on the cluster algebra gives rise to a well-defined
$\F^\star$-action on $\X$. The corresponding flow is called a {\it
toric flow}.

\proclaim{Lemma 2.3} Let $Z$ denote the transformation matrix
at a vertex $v$ of $\Bbb T_m$, and let $w$ be an arbitrary integer weight.
The local toric action at $v$ defined by $\{f_{1},\dots,f_{n}\}\mapsto
\{f_{1}t^{w_{1}},\dots,f_{n}t^{w_{n}}\}$ can be extended to a toric
action on $\cC$ if and only if $w^T$ belongs to the right nullspace
$N_{\er}(Z)$.
Moreover, if such an extension exists, then it is unique.
\endproclaim

\demo{Proof} Given a monomial $g=\prod_j f_{j}^{p_j}$, we define its
(weighted) degree by $\deg g=\sum_j p_jw_{j}$. It is easy to see that
local toric actions are compatible if and only if all the monomials in each
relation defining the transition $X\to \bar X$ have the same degree.

Consider an edge $(v,\bar v)$ of $\Bbb T_m$ labeled by $i$. Among the
relations defining the transition $X\to \bar X$ there is  a three-term
relation $\bar f_{i}f_{i}=M^i+\bar M^i$. Hence, the compatibility
condition implies $\deg M^i=\deg \bar M^i$, which by (1.1) is
equivalent to
$$
\sum_{z_{ik}>0} z_{ik}w_{k}=\sum_{z_{ik}<0}(-z_{ik})w_{k}.
$$
The latter condition written for all the $m$ edges incident to $v$ gives
$Zw^T=0$. Therefore, condition $w^T\in N_{\er}(Z)$ is necessary for the
existence of global extension.

Let us find the weight $\bar w$  that defines the local toric action at
$\bar v$ compatible with the initial local toric action at $v$. First,
identities $\bar f_{j}=f_{j}$ for $j\ne i$ immediately give
$\bar w_{j}=w_{j}$ for $j\ne i$. Next, the three-term relation gives
$$
\bar w_{i}=\sum_{z_{ik}>0} z_{ik}w_{k}-w_{i},
$$
so the weight at $\bar v$ is defined uniquely. It remains to prove that
$w^T\in N_{\er}(Z)$ implies $\bar w^T\in N_{\er}(\bar Z)$.

Let $k\ne i$, then the $k$th entry of $\bar Z\bar w^T$ is equal to
$$\align
\sum_{j=1}^n \bar z_{kj}\bar w_j
&=\sum_{j\ne i} z_{kj}w_j +\frac12\sum_{j\ne i}
\left(|z_{ki}|z_{ij}+z_{ki}|z_{ij}|\right)w_j
-z_{ki}\left(\sum_{z_{il}>0} z_{il}w_{l}-w_{i}\right)\\
&=\sum_{j=1}^n z_{kj}w_j=0.
\endalign
$$

The $i$th entry of $\bar Z\bar w^T$ is equal to
$$
\sum_{j=1}^n \bar z_{ij}\bar w_j=-\sum_{j\ne i} z_{ij}w_j=0,
$$
since $\bar z_{ii}=0$.
Hence,  $\bar w^T\in N_{\er}(\bar Z)$.
\qed
\enddemo

\subheading{2.3. Symplectic leaves}
Evidently, $\X$ is foliated into a disjoint union of symplectic leaves of the
Poisson brackets $\omega$. We are interested only in generic leaves,
which means the following.

Fix some generators $q_1,\dots,q_s$ of the field of rational Casimir functions $\F_C$.
They define a map $Q\:\X\to \F^s$, $Q(x)=(q_1(x),\dots,q_s(x))$.
Let $\L$ be a symplectic leaf, and let $z=(z_1,\dots,z_s)=Q(\L)\in\F^s$.
We say that $\L$ is {\it generic\/} if there exist $s$ vector fields $u_i$ in
a neighborhood of $\L$ such that

a) at every point $x\in \L$, the vector $u_i(x)$ is transversal to the
surface $q_i(x)=z_i$, which means that $\nabla_{u_i}q_i(x)\ne0$;

b) the translation along $u_i$ for a sufficiently small time $t$ gives a
diffeomorphism between $\L$ and a close symplectic leaf $\L_t$.

Let us denote by $\X^0$ the open part of $\X$ given by the conditions
$f_i\ne 0$ for $i\in [m+1, n]$. It is easy to see that
$\X^0=\cup_{v\in\Bbb T_m}X_v^0$.

\proclaim{Lemma 2.4} $\X^0$ is foliated into a disjoint union of
generic symplectic leaves of the Poisson bracket $\omega$.
\endproclaim

\demo{Proof} Consider first the special case when the
Poisson structure on $\X$ is nondegenerate at the generic point,
i.e., its rank equals to the dimension of the manifold. Then we
show that every point of $\X^0$ is generic, i.e., the rank of the
Poisson bracket is maximal at each point. Note that for every
point $x\in \X^0$ there exists a cluster chart $X_v$ such that on $X_v$
one has $f_{v,j}\ne 0$ for $j\in [1,n]$.
Therefore the coordinates $\log f_{v,j}$ form a local
coordinate system on $X_v$, and the Poisson structure written in these
coordinates becomes a constant Poisson structure. If a constant
Poisson structure is nondegenerate, it is nondegenerate at each
point, which proves the statement. Moreover, note that the
complement $\X\setminus \X^0$ consists of degenerate symplectic
leaves of smaller dimension. Hence, if the Poisson structure is symplectic,
i.e., nondegenerate at a generic point, then
$\X^0$ is a union of generic symplectic leaves.

Assume now that the rank of the Poisson structure is $r<n$. There
exist $s=n-r$ Casimir functions that generate the field $\F_C$. By Lemma~2.2,
one can build these Casimir functions by choosing $s$ independent integer vectors $u$ in the left nullspace
$N_{\el}(\Omega^\tau)$ and by constructing the corresponding
monomials $\m^u$. Observe that if $u=\pmatrix u^1&u^2\endpmatrix\in
N_{\el}(\Omega^\tau)$ and
$u'=\pmatrix D^{-1}u^1 &u^2\endpmatrix$, then $(u')^T\in N_{\er}(Z)$.
Therefore, by Lemma~2.3, such a
$u'$ defines a toric flow on $\X$.
To accomplish the proof it is enough to show that the toric flow
corresponding to the vector $u'$ is
transversal to the level surface $\{y\in \X^0 : \m^u(y)=\m^u(x)\}$, and that
a small translation along the trajectory of this toric flow transforms one
symplectic leaf into another one.

We will first
show that if $x(t)$ is a trajectory of the toric flow
corresponding to $u'$ with the initial value $x(1)=x$
and the initial velocity vector $\nu=d x(t)/dt|_{t=1}=(u'_1f_1,\dots,u'_nf_n)$
then $d \m^u(\nu)/dt\ne 0$.
Indeed, by Lemma~2.2,
$$
d \m^u(\nu)/dt=\sum_{i=m+1}^n \alpha_i \frac{\prod_{j=m+1}^n
f_{j}^{\alpha_j}}{f_{i}} \cdot u_if_i
=\m^u \alpha^2(u^2)^T.
$$
Since $x\in \X^0$, one has $\m^u(x)\ne0$. To find $\alpha^2(u^2)^T$ recall
that by the proof of Lemma~2.2, $\alpha^2=u^1Z_2+u^2Z_4+u^2K'$, where $K'$ is
the submatrix of $K$ whose entries lie in the last $n-m$ rows and columns.
Clearly, $u^2Z_4(u^2)^T=0$, since $Z_4$ is skew-symmetric. Next,
$u^1Z_2(u^2)^T=u^1\Lambda^{-1}\Omega_2(u^2)^T=u^1\Lambda^{-1}\Omega_1^T(u^1)^T=
0$, since $u^2\Omega_2^T=u^1\Omega_1$ and $\Lambda^{-1}\Omega_1^T$
is skew-symmetric. Thus, $\alpha^2(u^2)^T=u^2K'(u^2)^T\ne0$, since $K'$ can be
chosen to be a diagonal matrix with positive elements on the diagonal,
see the proof of Lemma~1.1.

Consider now another basis vector $\bar u\in N_{\el}(\Omega^\tau)$ and the
corresponding Casimir function $\m^{\bar u}$. It is easy to see that
$d \m^{\bar u}(\nu)/dt=\m^{\bar u}\bar\alpha^2(u^2)^T$.
Note that the latter expression does not depend on
the point $x$, but  only on the value $\m^{\bar u}(x)$ and on the vectors
$\bar\alpha$ and $u$. Therefore the value of the derivative
${d\m^{\bar u}}/{dt}$ is the same for all points $x$ lying on the same
symplectic leaf, and the toric flow transforms one symplectic leaf into
another one.
\qed
\enddemo

In general, it is not true that $\X^0$ coincides with the union of all
``generic" symplectic leaves. A simple counterexample is provided by the
cluster algebra given by two clusters $\{f_1,f_2,f_3\}$ and
$\{\bar f_1,f_2,f_3\}$ subject to one relation:
$\bar f_{1} f_{1}=f_{2}^2 f_{3}^2+1$.
One can choose the Poisson bracket on $\X$ as follows: $\{f_{1},f_2\}=
f_{1}f_2$, $\{f_{1},f_3\}=f_{1}f_3$, $\{f_2,f_3\}=0$. Equivalently,
$\{\bar f_{1},f_2\}=-\bar f_{1}f_2$, $\{\bar f_{1},f_3\}=-\bar f_{1}f_3$,
$\{f_2,f_3\}=0$. Generic
symplectic leaves of this Poisson structure are described by the
equation $A f_2+B f_3=0$ where $(A:B)$ is a homogeneous coordinate
on $P^1$. All generic symplectic leaves form $P^1$. In particular,
two leaves $(1:0)$ and $(0:1)$ (correspondingly, subsets $f_2=0,
f_3\ne 0, f_{1}\bar f_{1}=1$ and $f_3=0,f_2\ne 0,
f_{1}\bar f_{1}=1$) are generic symplectic leaves
in $\X$. According to the definition of $\X^0$ these leaves are not
contained in $\X^0$.

We can describe $\X^0$ as the nonsingular
locus of the toric action. The main source of examples of
cluster algebras are coordinate rings of homogeneous
manifolds. Toric actions on such cluster algebras are induced by the
natural toric actions on these manifolds.

All this suggests that $\X^0$ is
a natural geometrical object  in the cluster algebra
theory, intrinsically related to the corresponding Poisson structure.

\subheading{2.4. Connected components of $\X^0$}
In what follows we assume that $\F=\R$. In this case, the first natural
question concerning the topology of $\X^0$ is to find
the number $\#(\X^0)$ of connected components of $\X^0$.
To answer this question we follow the approach developed
in~\cite{SSV1, SSV2, Z}.

Given a vertex $v$ of $\Bbb T_m$, we define an open subset
$S(X^0)=S(X_v^0)\subset \X^0$ by
$$
S(X_v^0)=X_v^0\cup\bigcup_{(v,\bar v)\in\Bbb T_m}X_{\bar v}^0,
$$
where $(v,\bar v)\in\Bbb T_m$ means that $(v,\bar v)$ is an edge of
$\Bbb T_m$.

Recall that $X^0\simeq(\R\setminus 0)^n$. We can decompose $X^0$ as follows.
Let $\Sigma$ be the set
of all possible sequences $(\sigma(1),\dots,\sigma(n))$ of $n$
signs $\sigma(i)=\pm 1$. For $\sigma\in \Sigma$ we define
$X^0(\sigma)$ as the octant $ \sigma(j) f_{j}>0$ for all $j\in [1,n]$.
Two octants $X^0(\sigma_1)$ and $X^0(\sigma_2)$ are called {\it
essentially connected\/} if the following two conditions are fulfilled:

1) there exists $i\in [1,n]$ such that $\sigma_1(i)\ne \sigma_2(i)$
and $\sigma_1(j)= \sigma_2(j)$ for $j\ne i$;

2) there exists $x^*\in S(X^0)$ that belongs to the intersection of the
closures of $X^0(\sigma_1)$ and $X^0(\sigma_2)$.

The second condition can be restated as follows:

2') there exists $x^*\in S(X^0)$
such that $f_i(x^*)=0$, $\hat f_i(x^*)\ne 0$, and $f_j(x^*)\ne 0$ for $j\ne i$,
 where $\hat f=f_{\hat v}$ and
$(v,\hat v)$ is the edge of $\Bbb T_m$ labeled by $i$.

\proclaim{Lemma~2.5} Let $(v,\bar v)$ be an edge of $\Bbb T_m$, $\sigma_1,
\sigma_2\in \Sigma$, and let the octants $X^0(\sigma_1)$ and
$X^0(\sigma_2)$ be essentially connected. Then $\bar X^0(\sigma_1)$ and
$\bar X^0(\sigma_2)$ are essentially connected as well.
\endproclaim

\demo{Proof} Assume that the edge $(v,\bar v)$ is labeled by $j$, and
consider first the case $\sigma_1(j)\ne \sigma_2(j)$. Then $\hat v$ in
condition 2' coincides with $\bar v$, and hence $f_j(x^*)\bar f_j(x^*)=
M^j(x^*)+\bar M^j(x^*)=0$. Since $M^j$ and $\bar M^j$ both do not contain
$f_j$, any point $x$ such that
$f_l(x)=f_l(x^*)$ for $l\ne j$, $f_j(x)\ne0$, $\bar f_j(x)=0$ belongs
to $S(\bar X^0)$, and hence $\bar X^0(\sigma_1)$ and
$\bar X^0(\sigma_2)$ are essentially connected.

Assume now that $\sigma_1(i)\ne \sigma_2(i)$ for some $i\ne j$. As before,
we get $M^i(x)+\widehat M^i(x)=0$ for any $x$ such that $f_l(x)=f_l(x^*)$
for $l\ne i$. Consider the edge $(\bar v,\bar{\bar v})$ of $\Bbb T_m$
labeled by $i$; by the above assumption, $\bar{\bar v}\ne v$. Without
loss of generality assume that $\bar f_j$ does not enter $\bar{\bar M}^i$.
Then by~\cite{FZ2} one has
$$
M^i+\widehat M^i=\widehat M^i\left(\frac{M^i}{\widehat M^i}
+1\right)=\widehat M^i\left.\left(\frac{\bar M^i}{\bar{\bar M}^i}+1
\right)\right|_{\bar f_j\gets \frac{M_0}{f_j}},
$$
where $M_0=M^j+\bar M^j|_{f_i=0}$. Therefore, for $x^{**}$ such that
$\bar f_l(x^{**})=f_l(x^*)$ if $l\ne i, j$,
$\bar f_j(x^{**})=M_0(x^*)/f_j(x^*)$, $\bar f_i(x^{**})=0$, one has
$\bar M^i(x^{**})+{\bar{\bar M}^i}(x^{**})=M^i(x^*)+\widehat M^i(x^*)=0$,
and hence one can choose $\bar{\bar f}_i(x^{**})\ne0$.
\qed
\enddemo

\proclaim{Corollary 2.6} If $X^0_v(\sigma_1)$ and
$X^0_v(\sigma_2)$ are essentially connected, then
$X^0_{v'}(\sigma_1)$ and $X^0_{v'}(\sigma_2)$ are essentially
connected for any vertex $v'$ of $\Bbb T_m$.
\endproclaim

\demo{Proof} Since the tree $\Bbb T_m$ is connected, one can pick up
the path connecting $v$ and $v'$. Then the corollary follows immediately
from Lemma~2.5. \qed
\enddemo

Let $\#_v$ denote the number of connected components in $S(X_v^0)$.

\proclaim{Theorem 2.7} The number $\#_v$ does not depend on $v$ and is equal
to $\#(\X^0)$.
\endproclaim

\demo{Proof} Indeed,  since $S(X_v^0)$ is dense
in $\X^0$, one has $\#(\X^0)\ls \#_v$. Conversely, assume that
there are points $x_1,x_2\in \X^0$ that are connected by a path
in $\X^0$. Therefore their small neighborhoods are also connected
in $\X^0$, since $\X^0$ is a topological manifold. Since $X^0_v$ is
dense in $\X^0$, one can pick $\sigma,\sigma'\in\Sigma$ such that
the intersection of the first of the above neighborhoods with
$X_v^0(\sigma)$
and of the second one with $X_v^0(\sigma')$ are nonempty. Thus,
$X_v^0(\sigma)$ and $X_v^0(\sigma')$ are connected in $\X^0$.
Hence, there exist a loop $\gamma$ in $\Bbb T_m$ with the initial
point $v$, a subset $v_1,\dots,v_p$  of
vertices of this loop,  and a sequence
$\sigma_1=\sigma,\sigma_2,\dots,\sigma_{p+1}=\sigma'\in\Sigma$ such that
$X_{v_l}^0(\sigma_l)$ is essentially
connected with $X_{v_l}^0(\sigma_{l+1})$ for all $l\in [1,p]$. Then by
Corollary~2.6 $X_{v}^0(\sigma_l)$ and $X_{v}^0(\sigma_{l+1})$
are essentially connected. Hence all $X^0_v(\sigma_l)$ are connected
 with each other in $S(X^0_v)$. In particular,
$X_v^0(\sigma)$ and $X_v^0(\sigma')$ are connected in $S(X^0_v)$.
This proves the assertion.
\qed
\enddemo

By virtue of Theorem~2.7, we write $\#$ instead of $\#_v$.

Let $\FF_2^{n}$ be an $n$-dimensional vector space  over $\FF_2$ with a
fixed basis $\{e_i\}$. Let $Z'\in \Z^{\widehat D}_{nn}(Z)$ be as in
Section~1.2, and let $\eta=\eta_v$ be a (skew-)symmetric bilinear
form on $\FF_2^{n}$, such that $\eta(e_i,e_j)=d_iz'_{ij}$. Define a
linear operator $\frak t_i\:\FF_2^{n}\to \FF_2^{n}$ by the formula
$\frak t_i(\xi)=\xi-\eta(\xi,e_i) e_i$, and let $\Gamma=\Gamma_v$ be the
group generated by $\frak t_i$, $i\in [1,m]$.

The following lemma is a minor modification of the result presented in
\cite{Z}.

\proclaim{Lemma~2.8} The number of connected components
$\#_v$ equals to the number of $\Gamma_v$-orbits
in $\FF_2^{n}$.
\endproclaim

For reader's convenience we will repeat the proof of this
statement here.

Let us identify $\FF_2^{n}$ with $\Sigma$ by the following rule:
a vector $\xi\in\FF_2^{n}$ corresponds to $\sigma\in\Sigma$ such that
$\sigma(i)=(-1)^{\xi_i}$. Abusing notation we will also write
$X^0(\xi)$ instead of $X^0(\sigma)$.

\proclaim{Lemma 2.9} Let $\xi$ and $\xi'$ be two
distinct vectors in $\FF_2^{n}$. Then the closures of
$X^0(\xi)$ and $X^0(\xi')$ intersect in $S(X^0)$
if and only if $\xi' = \frak t_i (\xi)$ for
some $i\in [1,m]$.
\endproclaim

\demo{Proof} Suppose $x \in S(X^0)$ belongs to the intersection
of the closures of $X^0(\xi)$ and $X^0(\xi')$.
Then $f_l (x) = 0$ whenever $\xi_l \neq \xi'_l$.
From the definition of $S(X^0)$ we see that there is a unique
$i$ such that $\xi_i \neq \xi'_i$; evidently, $i\in [1,m]$.
Furthermore, if $(v,\bar v)$ is an edge of $\Bbb T_m$ labeled by
$i$, then  $\bar f_{i}(x) \neq 0$. Since
any neighborhood of $x$ intersects both $X^0(\xi)$ and
$X^0(\xi')$, it follows that monomials $M^i$ and
$\bar M^i$ on the right hand side of the $3$-term relation
$f_{i}\bar f_{i}=M^i+\bar M^i$ must have opposite signs
at $x$. Therefore
$$\xi_i -
\xi'_i = 1 = \sum_{j=1}^n d_iz_{ij} \xi_j.
$$
Comparing this with the definition of the transvection $\frak t_i$,
we conclude that $\xi' = \frak t_i (\xi)$, as claimed.

Conversely, suppose $\xi' = \frak t_i (\xi) \neq
\xi$, then $\sum_{j=1}^n d_iz_{ij} \xi_j=1$.
Therefore, there exists a point $x \in S(X^0)$ such that
$(-1)^{\xi_l} f_l (x) > 0$ for all $l \neq i$, and the right hand
side of three $3$-term relation
vanishes at $x$. Then any neighborhood of
$x$ contains points where the signs of all $f_l$ for $l \neq i$
remain the same, while the right hand side of the $3$-term relation
is positive (or negative). Thus, $x$ belongs to the intersection
of the closures of $X^0(\xi)$ and $X^0(\xi')$, and we are done.
\qed
\enddemo

Now we are ready to complete the proof of
Lemma~2.8.
Let $\Xi$ be a $\Gamma$-orbit in $\FF_2^{n}$, and
let $X^0(\Xi)\subset S(X^0)$ be the union of the closures of $X^0(\xi)$
over all $\xi\in \Xi$. Each
$X^0(\xi)$ is a copy of $\R_{> 0}^{n}$, and is thus connected.
Using the ``if" part of Lemma~2.9,
we conclude that $X^0(\Xi)$ is connected (since the closure of a
connected set and the union of two non-disjoint connected sets are
connected as well). On the other hand, by the ``only if" part of
the same lemma, all the sets $X^0(\Xi)$ are pairwise disjoint.
Thus, they are the connected components of $S(X^0)$, and
we are done.
\qed
\medskip

Theorem~2.7 and Lemma~2.8 imply the following theorem.

\proclaim{Theorem~2.10} The number of connected components
$\#(\X^0)$ equals the number of $\Gamma$-orbits in
$\FF_2^{n}$.
\endproclaim

\heading 3. Poisson and cluster algebra structures on Grassmannians \endheading

\subheading{3.1. Sklyanin bracket and factorization parameters}
As we have seen from the discussion above, one can associate
a cluster algebra with any algebraic Poisson manifold equipped with a
system of rational log-canonical coordinates.  A rich collection
of non-trivial examples of this sort is provided by the theory of Poisson-Lie
groups
and Poisson homogeneous spaces. This collection includes real Grassmannians,
which will serve as our main example.

Recall (see e.g.~\cite{ReST}) that a Lie group $G$ equipped with a Poisson
structure is called
a Poisson-Lie group if the multiplication map $G\times G \to G$ is Poisson.
First, we review  the definition of the Sklyanin Poisson bracket and the
standard Poisson-Lie structure on a semi-simple Lie group $G$.

Let $B_+$ and $B_-$ be two $\R$-split
opposite Borel subgroups, $N$ and $N_-$ be their unipotent radicals,
$H = B_+ \cap B_-$ be an $\R$-split maximal torus of $G$, and $W =
\Norm_G(H)/H$ be the Weyl group of $G$. For every $x\in N_- H N_+$,
we write its unique Gauss factorization as $x=x_- x_0 x_+$
and define $x_{\ls 0} =x_- x_0$, $x_{\gs 0} = x_0 x_+$.

As usual, let $\g$, $\h$, $\b_\pm$, $\n_\pm$ be Lie algebras that correspond
to  $G$, $H$, $B_\pm$, $N_\pm$.
We denote by $\Phi$ the root system of $\g$, by
$\Phi^+$ (resp. $\Phi^-$) the set of all positive (resp. negative)
roots, and by $\alpha_1,\dots, \alpha_r$
simple positive roots. To each $\alpha_i$ there corresponds
the elementary reflection $s_i\in W$. A
{\it reduced decomposition\/} for $w \in W$ is a factorization of $w$
into a product $w = s_{i_1} \cdots s_{i_l}$ of simple reflections,
where $l$ is the smallest length of such a factorization. Then an integral
vector
$\ii=(i_1,\ldots,i_l)$ is called a {\it reduced word\/} that corresponds
to $w$, while $l$ is called the
{\it length\/} of $w$ and is denoted by $l(w)$.

We fix a Chevalley basis
$\{ e_\alpha, \alpha\in\Phi; \ h_i, i\in [1,r]\}$ in $\g$.
The Killing form on $\g$ will be denoted by $\langle\ ,\ \rangle$.

Let $R\in \End(\g)$ be a skew-symmetric map. The Sklyanin
bracket is defined by
$$
\omega(f_1, f_2) (x)=
\langle R (\nabla_{\text{\eurs r}} f_1 (x)), \nabla_{\text{\eurs r}} f_2 (x) \rangle -
\langle R (\nabla_{\text{\eurs l}} f_1 (x)) , \nabla_{\text{\eurs l}} f_2 (x) \rangle.
\tag 3.1
$$
Here  the right and left gradients $\nabla_{\text{\eurs r}}$, $\nabla_{\text{\eurs l}}$
are defined
by
$$
\langle \nabla_{\text{\eurs r}} f(x),\xi\rangle=
\left.\frac{d}{dt}f(x\exp(t\xi))\right|_{t=0},
\qquad\langle\nabla_{\text{\eurs l}} f(x),\xi\rangle=\left.\frac{d}{dt} f(\exp(t\xi)x)
\right|_{t=0}
$$
for $x\in G$, $\xi \in \g$.

If $R$ satisfies the modified classical Yang-Baxter equation (MCYBE)
$$
[R \xi, R \eta]- R \left ( [R \xi, \eta] + [ \xi, R \eta] \right )
= - \alpha [\xi,\eta], \quad \xi,\eta \in \g,
$$
then the bracket (3.1) satisfies the Jacobi identity and equips $G$
with a structure of the Poisson-Lie group.

For any $\xi \in \g$ we write its unique decomposition
into a sum of elements in $\n_-, \h$ and $\n_+$ as
$$
\xi=\xi_- + \xi_0 + \xi_+.
$$
The standard Poisson-Lie structure on $G$ denoted by  $\omega_G$ corresponds
to a particular solution of MCYBE:
$$
R(\xi)= \xi_+ - \xi_-
\tag 3.2
$$
and can be characterized as follows \cite{HKKR, R}.

For $t\in \C$, denote $x^+_i(t)=\exp(t e_{\alpha_i})$, $x^-_i(t)=
\exp(t e_{-\alpha_i})$.
For every $i\in [1,r]$ one defines a group homomorphism
(canonical inclusion) $\rho_i$
from $SL_2$ to $G$ by
$$
\rho_i\pmatrix 1  & t \\ 0 & 1 \endpmatrix=x^+_i(t),\qquad
\rho_i\pmatrix 1  & 0 \\ t & 1 \endpmatrix=x^-_i(t).
$$
Let $D_i$ be the length of the $i$th simple root $\alpha_i$. Then $\omega_G$
 is the unique
Poisson-Lie structure on $G$ such that each map
$$
\rho_i\: \left(SL_2, D_i \omega_{SL_2} \right) \to (G, \omega_G)
\tag 3.3
$$
is Poisson.

Symplectic leaves of $\omega_G$ have been studied in \cite{HL}.
Their explicit description
was obtained in more recent works \cite{HKKR, R, Y, KoZ}, where the key role
is played by {\it double Bruhat cells\/} that were comprehensively studied
in \cite{FZ1}.
We shall now recall the definition of a double Bruhat cell and review the
results
of \cite{HKKR, R, KoZ} that show that each such cell is equipped with a family
of log-canonical coordinates of the kind described in the first part of this
paper.

Recall that the group $G$ has two {\it Bruhat decompositions}, with respect
to $B_+$ and $B_-$:
$$
G = \bigcup_{u \in W} B_+ u B_+ = \bigcup_{v \in W} B_- v B_-.
$$
The double Bruhat cells $G^{u,v}$ are
defined by $G^{u,v} = B_+ u B_+  \cap B_- v B_-$. According to \cite{FZ1},
the variety $G^{u, v}$ is biregularly isomorphic to a Zariski open
subset of $\C^{r+l(u)+l(v)}$. Furthermore, if for every pair $\jj$, $\kk$
of reduced words for $u$
and $v$
one defines a word $\ii=(i_\nu)_{\nu=1}^{m}$ with
$m=l(u)+l(v)$ as an
arbitrary shuffle of the words $-\jj$ and  $\kk$, then the map
$x_{\ii}\: H\times \C^{m}\to G^{u, v}$ defined by
$$
x_{\ii}(h,t)= h
\prod_{\nu=1}^{m}
x^{\s(i_\nu)}_{i_\nu}(t_\nu),
\tag 3.4
$$
where $h\in H$ and $t=(t_{1},\ldots,t_{m})\in \C^m$,
restricts to a biregular isomorphism between $H\times \C_{\ne 0}^{m}$
and a Zariski open subset of $G^{u, v}$. Let us further factor $h$ in
(3.4)
into $h= \rho_1(\diag(a_1,a_1^{-1}))\cdots \rho_r(\diag(a_r,a_r^{-1}))$.
Parameters $a_1,\ldots,a_r; t_{1},\ldots,t_{m}$ are called {\it factorization
parameters}. Explicit
formulae for the inverse of the map (3.4) were found in \cite{FZ1}.

The relevance of double Bruhat cells and factorization parameters in the
context of the standard Poisson-Lie structure was observed in
\cite{HKKR, R}. It turns out
that (i) for every $u, v$, $G^{u,v}$ is a Poisson submanifold of
$(G, \omega_G)$, and
(ii) for every $\ii$, factorization parameters form a family of log-canonical
coordinates
on a Zariski open subset of $G^{u, v}$. Both assertions can be verified via
the following
construction \cite{HKKR, R}. For $i\in [1,r]$ define subgroups
$B_i^\pm$ of $B_\pm$
as images of resp. upper and lower triangular subgroups of $SL_2$ under the
homomorphism $\rho_i$:
$$
B_i^+= \rho_i \left(\left\{ \pmatrix d  & c \\ 0 & d^{-1} \endpmatrix\right \}
\right ),\qquad
B_i^-= \rho_i \left(\left\{ \pmatrix d  & 0 \\ c & d^{-1} \endpmatrix\right \}
\right ).
\tag 3.5
$$
We can view parameters $c$, $d\ne 0$  as coordinates on $B_i^+$
(resp. $B_i^-$) and
denote the corresponding elements of $B_i^\pm$ by $b_i^\pm(c,d)$.
In view of the Poisson property (3.3) of the map $\rho_i$, one obtains
for both $B_i^+$ and $B_i^-$
$$
\omega(d, c)= D_idc.
\tag 3.6
$$
For any $u, v\in W$ and the corresponding fixed $\ii$, denote
$$
B_\ii = \times_{\nu=1}^{m} B_{i_{\nu}}^{\s(i_{\nu})}
\tag 3.7
$$
and define the multiplication map $y_\ii \: B_\ii \to G $ by
$$
y_{\ii}(b_\ii)=
\prod_{\nu=1}^{m}b^{\s(i_\nu)}_{i_\nu}
(c_\nu, d_\nu),
\tag 3.8
$$
where
$b_\ii =
\left (b_{i_1}^{\s(i_1)}(c_1,\ d_1), \ldots, b_{i_{m}}^{\s(i_{m})}
(c_{m},\ d_{m})\right)\in B_{\ii}$.
This map is clearly Poisson. Moreover, comparing (3.8) and (3.4),
one sees that the image of the restriction of $y_\ii$ to the set where all
$b_\alpha\ne 0$
is a Zariski open subset of $G^{u,v}$ that
coincides with the image of restriction of $x_\ii$ to  $H\times\C_{\ne 0}^{m}$.
It is not hard to see that factorization parameters $a_1,\ldots,a_r; t_{1},\ldots,t_{m}$
are monomial functions in parameters $c_\nu, d_\nu$ and since by (3.6)
the latter are log-canonical, the former are log-canonical as well.

The only drawback of the log-canonical coordinate system
$a_1,\ldots,a_r; t_{1},\ldots,t_{m}$ is that factorization
parameters are rational functions on $G^{u, v}$. However, one of
the main theorems of \cite{FZ1} gives an explicit expression for
these parameters as an  invertible monomial transformation of a
family of functions {\it regular\/} on $G^{u, v}$, the so-called {\it
twisted generalized minors}. For every $\ii$,
log-canonical Poisson brackets of functions in this family w.r.t.
the standard Poisson-Lie structure were computed explicitly
in a recent paper \cite{KoZ} and used to refine a
description of symplectic leaves of the Poisson-Lie group $G$. The
definition and properties of generalized minors of an element
$x\in G$ can be found in \cite{FZ1} and will not be reproduced
here (in the
$SL_n$ case generalized minors are the minors of a matrix $x$).
Twisted generalized minors are generalized minors of the {\it
twist map\/} $ G^{u,v}\ni x \mapsto x'\in G^{u^{-1},v^{-1}}$ defined
in \cite{FZ1} as follows. For any $w = s_{i_1} \cdots s_{i_l}\in W$,
pick a representative $\widehat{w}$ of $w^{-1}$ in $G$ given by
$$
\widehat{w} =
\prod_{\nu=1}^{l} \rho_{i_{l+1-\nu}}
\pmatrix 0  & -1 \\ 1 & 0 \endpmatrix
$$
and define
$$
x'=\left ( (\widehat{u} x)^{-1}_- \widehat{u} x\widehat{v}
(x\widehat{v})^{-1}_+\right )^\theta,
\tag 3.9
$$
where $\theta$ is the involutive automorphism of $G$ uniquely determined
by
$$
x^\pm_i(t)^\theta = x^\mp_i(t), \qquad h^\theta= h^{-1},\quad h\in H
$$
and $a_-^{-1}$ stands for $(a_-)^{-1}$.

We conclude this section with the result that shows that a
construction  of \cite{KoZ} can be used to build another system
of log-canonical regular coordinates on $G^{u, v}$ from
generalized minors of $x$ rather than its twist $x'$.

\proclaim{Theorem 3.1}
The twist map {\rm (3.9)} is an anti-Poisson map with respect to
the Sklyanin bracket {\rm (3.1), (3.2)}.
\endproclaim

\demo{Proof} First, observe that an automorphism $\theta$
itself is an anti-Poisson map from $G$ to $G$. This can be easily
seen from the $\theta$'s action on subgroups $B^\pm_i$ defined in
(3.5). Thus, we only need to check that the map $x\mapsto\chi_{u,v}(x)$,
where $\chi_{u,v}(x)$ is defined by the expression inside the outer brackets
in (3.9), is Poisson. It suffices to check this property for elements of a
Zariski open subset of $G^{u,v}$. Let $G_0^{u,v}=G^{u,v} \cap N_- H N_+$.
For $x\in G_0^{u,v}$,
$$\align
\chi_{u, v} (x)&=
(\widehat{u} x)^{-1}_- (\widehat{u} x)x^{-1}(x\widehat{v})
(x\widehat{v})^{-1}_+=
(\widehat{u} x)_{\gs 0}\ x^{-1}(x\widehat{v})_{\ls 0}\\
&=(\widehat{u} x_{\ls 0})_{\gs 0}x_+x^{-1}x_- (x_{\gs 0}\widehat{v})_{\ls 0}
=(\widehat{u} x_{\ls 0})_{\gs 0}\ x_0^{-1}
(x_{\gs 0}\widehat{v})_{\ls 0}.
\endalign
$$
Furthermore, fix reduced words $\jj= (j_\nu)_{\nu=1}^{l(u)}$,
$\kk=(k_\gamma)_{\gamma=1}^{l(v)}$
for $u$ and $v$ and assume that $x$ belongs
to the image of the multiplication map
$y_\ii\: B_\ii = B_{-\jj}\ \times B_\kk  \to G^{u,v}$
defined as in (3.8), (3.7),
where $\ii=(-\jj,\kk)$. Any element $b$ in the preimage of $x$ under the
map $y_\ii$ can be written as $b=(b_1,b_2)$, where $b_1\in B_{-\jj}$,
$b_2\in B_{\kk}$, and hence $x=y_\ii(b)=y_{-\jj}(b_1)y_\kk(b_2)=
x_-d_1d_2x_+$, where $d_1=(y_{-\jj}(b_1))_0$, $d_2=(y_\kk(b_2))_0$.
Therefore,
$$\align
\chi_{u,v}(y_\ii(b)) &= (\widehat{u} x_{\ls 0}d_2^{-1})_{\gs 0}
(d_1^{-1}x_{\gs 0}\widehat{v})_{\ls 0}=
(\widehat{u} y_{-\jj}(b_1))_{\gs 0}(y_{\kk}(b_2)\widehat{v})_{\ls 0}\\
&=\chi_{u,\id } (y_{-\jj}(b_1))\chi_{\id,v}(y_\kk(b_2) ).
\endalign
$$
Denote $\tilde \jj = (j_{l(u)-\nu+1})_{\nu=1}^{l(u)}$,
$\tilde \kk = (k_{l(v)-\gamma+1})_{\gamma=1}^{l(v)}$. It follows from the above
equality that
in order to prove that the map $x\mapsto \chi_{u,v}(x)$ is Poisson it is enough
to show that maps $\chi_{u,\id}\: G^{u,\id} \to G^{\id,u^{-1}}$ and
$\chi_{\id,v}\:G^{\id,v} \to G^{v^{-1},\id }$ are Poisson. Both maps can be
treated in the same fashion, so we shall concentrate on the second one.

Let $x\in G^{\id,v}=B_+ \cap B_- v B_- $, put $q=i_{l(v)}$ and write $v$ as
$v= w s_{q}$. Assume further that $x$ is in the image of the
multiplication map $G^{\id,w}\times B_q^+ \ni (x_1, b^+_q(c,d))\mapsto
x_1 b^+_q(c,d)\in G^{\id,v}$. (The map is Poisson and its image
is Zariski open in $G^{\id,v}$.) We shall show that
$\chi_{\id,v}(x) = (\tilde x\widehat{w})_{\ls 0}$, where
$\tilde x$ belongs to $G^{s_q,w}$ and the map $x\mapsto \tilde x$ is
Poisson. Then the statement will follow by induction on the length of $v$.

Observe that
$$\align
x \widehat{v} &=
x_1 \left ( b^+_q(c,d) \widehat{s_q}\right )\widehat{w}=
x_1 \left ( b^-_q(d^{-1},\ c)  b^+_q(-d, 1)\right )\widehat{w}\\
&=\left (x_1  b^-_q(d^{-1},\ c)\right ) \widehat{w}
\left ( (\widehat{w})^{-1} b^+_q(-d, 1)
\widehat{w}\right ).
\endalign
$$
So we can define $\tilde x = x_1  b^-_q(d^{-1},c)$ and note that
$ l(s_q w^{-1})
= l(w^{-1}) +1$ implies
$$
\widehat{w}^{-1} b^+_q(-d, 1)\widehat{w} \in N_+.
$$
Thus, $\chi_{\id,v}(x) = (\tilde x\widehat{w})_{\ls 0}$, and
since, by
(3.6), the map $b^+_q(c,d) \mapsto b^-_q(d^{-1},c)$ is Poisson, so is the map
$x \mapsto \tilde x$. This completes the proof.
\qed
\enddemo

As we have seen from the discussion above, the real part of any double
Bruhat cell $G^{u,v}$
can be equipped (in more than one way) with regular log-canonical coordinates
that serve
as a starting point for a construction of a cluster algebra. Results of
\cite{FZ1, FZ2, HKKR, R, Y} show that both Poisson and cluster algebra
structure on $G^{u,v}$ are determined by the Lie group structure of $G$ and
are compatible
in the way described in the first part of this paper.
We now turn to our main example, real Grassmannians, to show how a structure
of a Poisson
homogeneous space can also lead to a construction of a cluster algebra.

\subheading{3.2. Poisson structure and log-canonical coordinates}
Let $\P$ be a Lie subgroup of a Poisson-Lie group $(G, \omega_G)$.
A Poisson structure on
the homogeneous space $\P \backslash G$ is called {\it Poisson homogeneous\/}
if the action map
$\P \backslash G \times G \to \P \backslash G$ is Poisson. Conditions on
$\P$ for the Sklyanin
bracket (3.1) to descend to a Poisson homogeneous structure on
$\P \backslash G$
are conveniently formulated it terms of the Manin triple that corresponds to
$(G, \omega_G)$
and can be found, e.g., in \cite{ReST}. In particular, these conditions are
satisfied for $G=SL_n$  and
$$
\P=\P_k = \left \{ \pmatrix A & 0\\ B & C \endpmatrix \:  A\in
SL_k,  C\in SL_{n-k} \right \}.
$$
The resulting homogeneous space is the Grassmannian $G(k,n)$.

In what follows, we will need an explicit expression of the Poisson homogeneous
brackets on $G(k,n)$. First, consider the Sklyanin bracket (3.1), (3.2)
on $SL_n$. The form $\langle\ ,\ \rangle$ now coincides with the trace form:
$$
\langle A, B \rangle =\Trace AB.
$$
The Sklyanin bracket can be extended from $SL_n$ to the associative algebra
$\Mat_n$ of $n\times n$ matrices;  it is given there by
$$
\omega(f_1, f_2)(X)=
\langle R (\nabla f_1 (X) X ), \nabla f_2 (X) X \rangle -
\langle R (X \nabla f_1 (X)) , X \nabla f_2 (X) \rangle,
\tag 3.10
$$
where the gradient $\nabla$ is defined w.~r.~t.~the trace form.
In terms of matrix elements $x_{ij}$, $i,j\in [1, n]$, of a matrix
$X\in \Mat_n$, (3.10) looks as follows:
$$
\omega(x_{ij}, x_{\alpha\beta})= (\s (\alpha - i)+ \s (\beta -j))
x_{i\beta} x_{\alpha j}.
$$

If $X\in SL_n$ admits a factorization into block-triangular matrices
$$
X= \pmatrix X_1 & 0\\ Y' & X_2 \endpmatrix
\pmatrix \one_k &  Y \\ 0 & \one_{n-k} \endpmatrix = V U,
$$
then $Y$ represents an element of the cell $G_0(k,n)$ in $G(k,n)$
characterized
by non-vanishing of the Pl\"ucker coordinate $\pi_{[1,k]}$.

Relations between the Pl\"ucker coordinates $\pi_I$, $I=\{i_1,\ldots,i_k\:
1\ls i_1 < \cdots < i_k \ls n\}$, and minors
$Y_{\alpha_1 \ldots \alpha_l}^{\beta_1 \ldots \beta_l}=
\det(y_{\alpha_i,\beta_j})_{i,j=1}^l$
of $Y$ are given  by
$$
Y_{\alpha_1 \ldots \alpha_l}^{\beta_1 \ldots \beta_l} =
(-1)^{k l -l(l-1)/2 - (\alpha_1 + \cdots + \alpha_l)}
\frac{\pi_{([1,k]\setminus\{\alpha_1 \ldots \alpha_l\})
\cup \{\beta_1+k \ldots \beta_l+k\}}}{\pi_{[1,k]}}.
$$
Note that, if the row index set $\{\alpha_1 \ldots \alpha_l\}$ in the
above formula
is contiguous then the sign in the right hand side can be expressed
as $(-1)^{(k-\alpha_l) l}$.

It is not hard to see that the variation of $Y=Y(X)$ is given by
$$
\delta Y = \pmatrix \one_k & 0 \endpmatrix
 V^{-1} \delta X U^{-1}
\pmatrix 0 \\ \one_{n-k} \endpmatrix,
$$
and, therefore,
$$
\nabla (f\circ Y) = U^{-1} \pmatrix 0 \\ \one_{n-k} \endpmatrix
\nabla f \pmatrix \one_k & 0 \endpmatrix V^{-1}
$$
and
$$\gather
\nabla (f\circ Y) X = \Ad_{U^{-1}} \left(\pmatrix 0 \\ \one_{n-k} \endpmatrix
\nabla f \pmatrix \one_k & 0 \endpmatrix \right )  =
\pmatrix -Y \\ \one_{n-k}\endpmatrix \nabla f \pmatrix\one_k & Y \endpmatrix,\\
 X \nabla (f\circ Y) = \Ad_{V} \left ( \pmatrix 0 \\ \one_{n-k} \endpmatrix
\nabla f \pmatrix \one_k & 0 \endpmatrix \right ) \in \n_-.
\endgather
$$
Thus we obtain from (3.10)
$$\gather
\omega(f_1\circ Y,f_2\circ Y) =
\left \langle  \pmatrix  - R (Y \nabla f_1) & -Y \nabla f_1 Y\\
- \nabla f_1 &  R (\nabla f_1 Y) \endpmatrix,
\pmatrix - Y \nabla f_2 & -Y \nabla f_2 Y\\
\nabla f_2 &  \nabla f_2 Y \endpmatrix
\right \rangle\\
=\langle R (\nabla f_1  Y ), \nabla f_2 Y\rangle + \langle R (Y
\nabla f_1) , Y \nabla f_2)] \rangle = \omega_{G(n,k)}(f_1,f_2)\circ Y.
\endgather
$$

In terms of matrix elements $y_{ij}$, this formula  looks as follows:
$$
 \omega(y_{ij}, y_{\alpha\beta})= (\s (\alpha - i)- \s (\beta -j))
y_{i\beta} y_{\alpha j}. \tag 3.11
$$

Next, we will introduce new coordinates on $G(k,n)$, log-canonical
w.~r.~t.~$\omega$. This
will require some preparation.

Let $I=\{i_1,\ldots, i_r\}, J=\{j_1,\ldots, j_r\}$ be ordered
multi-indices. We denote by $I(i_p\to \alpha)$ the result of
replacing $i_p$ with $\alpha$ in $I$, by $I\setminus i_p$ the multi-index
obtained by deleting $i_p$ from $I$ and by $(\alpha I)$ the
multi-index $I=\{\alpha, i_1,\ldots, i_r\}$. For a matrix $X$, we denote
$X_I^J=\det X(I;J)=\det (x_{i_p j_q})_{p,q=1}^r$.
Then  the Laplace expansion formula implies
$$
\sum_{p=1}^r x_{i_p \beta} X_{I(i_p\to \alpha)}^J = x_{\alpha
\beta} X_I^J - X_{(\alpha I)}^{(\beta J)}.
\tag 3.12
$$

We will say that $\alpha < I$ (resp. $\alpha > I$), if
 $\alpha$ is less than the minimal index in $I$ (resp., the maximal index in $I$ is
less than $\alpha$).
We define
$\s (\alpha - I) = -\s  (I-\alpha)$ to be $-1, 0$ or $1$, if
$\alpha < I, \alpha \in I $ or $\alpha > I$, resp. Otherwise, $\s (\alpha - I)$
is not defined.

\proclaim{Lemma 3.2}
If $\s (\alpha -I)$ and $\s (\beta -J)$ are defined
and
$$
|\s (\alpha -I)- \s (\beta -J)| \ls 1, \tag 3.13
$$
then
$$
\omega(y_{\alpha \beta}, Y_I^J) =
- \left (\s (\alpha -I)- \s (\beta -J)\right )y_{\alpha \beta} Y_I^J.
\tag 3.14
$$
\endproclaim

\demo{Proof}
It is evident from (3.11) that
$\omega(y_{\alpha \beta}, Y_I^J) = 0$ if $\alpha < I$, $\beta < J$,
or $\alpha > I$, $\beta > J$,
so in these cases (3.14) holds true.

In general, one obtains from (3.11)
$$\align
\omega(y_{\alpha \beta}, Y_I^J) &=\sum_{p,q=1}^r (-1)^{p+q}\{y_{\alpha \beta},
 y_{i_p}^{j_q} \}
Y_{I\setminus i_p}^{J\setminus j_q}\\
&=\sum_{p=1}^r \s (i_p - \alpha) y_{i_p \beta}
\sum_{q=1}^r (-1)^{p+q}y_{\alpha j_q}
Y_{I\setminus i_p}^{J\setminus j_q}\\
&\quad - \sum_{q=1}^r \s (j_q - \beta) y_{\alpha j_q}
\sum_{p=1}^r (-1)^{p+q}y_{i_p \beta}
Y_{I\setminus i_p}^{J\setminus j_q}\\
&=\sum_{p=1}^r \s (i_p - \alpha) y_{i_p \beta} Y_{I(i_p\to\alpha)}^J
-\sum_{q=1}^r \s (j_q - \beta) y_{\alpha j_q} Y_I^{J(j_q\to\beta)}\ .
\endalign
$$

Assume that $\beta \in J$. Then, in the second sum above,
$Y_I^{J(j_q\to\beta)}$ can be nonzero only if $\s (j_q - \beta)=0$, which
implies that the sum is zero. If in addition $\alpha \in I$,
the first sum is equal to zero as well, and thus
$\omega(y_{\alpha \beta}, Y_I^J ) = 0$ if $\alpha \in I, \beta \in J$, which
is consistent with  (3.14). Otherwise, by our assumptions,
$\s (i_p - \alpha)$ does not depend on $p$ and is equal to $\s (I - \alpha)$.
In this case, (3.12) implies
$$
\omega(y_{\alpha \beta}, Y_I^J) = \s (I - \alpha) (y_{\alpha
\beta} Y_I^J - Y_{(\alpha I)}^{(\beta J)} )= \s (I - \alpha) y_{\alpha
\beta} Y_I^J.
$$
This agrees with (3.14). The remaining case $\alpha \in I, \s (J - \beta)=\pm 1$
can be treated in the same way.
\qed
\enddemo

The following Corollary drops out immediately from Lemma~3.2
together with the Leibnitz rule for Poisson brackets.

\proclaim{Corollary 3.3} Let $A=\{ \alpha_1,\ldots,\alpha_l\} , B=\{\beta_1,\ldots,\beta_l\}$
be such that for every pair $(\alpha_p, \beta_q)$, $p,q\in [1,l]$, condition
{\rm(3.13)} is satisfied. Then
$$
\omega(Y_A^B, Y_I^J) = - \left ( \sum_{p=1}^l ( \s (\alpha_p -I)- \s (\beta_p -J))\right )
Y_A^B Y_I^J.
\tag 3.15
$$
\endproclaim

For every $(i,j)$-entry of $Y$ define
$l(i,j)=\min (i-1, n-k-j)$
and
$$
F_{ij}= Y_{i-l(i,j),\ldots,i}^{j,\ldots, j+l(i,j)}.
\tag3.16
$$
It is easy to see that the change of coordinates $(y_{ij})\mapsto(F_{ij})$ is a
birational transformation.

\proclaim{Proposition 3.4}
Put
$$
t_{ij}= \frac{ F_{ij}}{ F_{i-1,j+1}}.
\tag3.17
$$
Then
$$
\omega_{G(n,k)}(\ln t_{ij} ,\ln t_{\alpha \beta}) =
\s (j-\beta)\delta_{i\alpha}- \s (i-\alpha)\delta_{j\beta}.
\tag 3.18
$$
\endproclaim

\demo{Proof}
First, we will show that coordinates $F_{ij}$ are log-canonical.
For this one needs to check that conditions of Corollary~3.2
are satisfied for every pair $F_{\alpha\beta}, F_{ij}$.
One has the following seven cases to consider.

1) $\alpha \ls i$,  $\beta\ls j$,  $i+j\ls n-k +1$;

2) $\alpha \ls i$,  $\beta > j$,  $\max(\alpha+\beta, i+j)\ls n-k +1$;

3) $\alpha \ls i$,  $\beta\ls j$,  $\alpha+\beta > n-k +1$;

4) $\alpha \ls i$,  $\beta > j$,  $\min(\alpha+\beta, i+j) > n-k +1$;

5) $\alpha \ls i$,  $\beta\ls j$,  $\alpha+\beta \ls n-k +1 < i+j$;

6) $\alpha \ls i$,  $\beta > j$,  $\alpha+\beta \ls n-k +1 < i+j$;

7) $\alpha \ls i$,  $\beta > j$,   $i+j \ls n-k +1 < \alpha+\beta$.

Direct inspection shows that choosing $Y_A^B=F_{ij}$, $Y_I^J=F_{\alpha\beta}$ in  case 3
and $Y_A^B=F_{\alpha\beta}$, $Y_I^J=F_{ij}$ in all the remaining cases
ensures that conditions of Corollary~3.3 hold true.
Moreover, one can use (3.15) to compute
$\omega_{G(n,k)}(\ln F_{\alpha \beta}, \ln F_{ij})$ in each of these cases;
the answers are

1) $- \min(\alpha, j - \beta)$;

2) $\max(\alpha+\beta, i+j) - \max(\beta, i+j)$;

3) $i - \max (\alpha , i+j+k -n-1))$;

4) $\min(\alpha, i+j+k-n-1) - \min(\alpha+\beta+k-n-1, i+j+k-n-1)$;

5)  $n+1 -j - k + \max(\alpha, i+j+k-n-1) - \max(\alpha+\beta, i)$;

6) $-(n+1 - k - \max(\beta, i+j))$;

7)  $-\min(\alpha, i+j+k-n-1)$.

Note that formulae above remain valid if we replace all strict inequalities
used to describe cases 1 through 7 by non-strict ones.
Now (3.18) can be derived from the formulae above via the case by case
verification, simplified by noticing that due to the previous remark, for any
$(i,j)$ and $(\alpha,\beta)$, all the four quadruples $(i,j),(\alpha,\beta)$;
$(i-1,j+1),(\alpha,\beta)$; $(i,j),(\alpha-1,\beta+1)$ and $(i-1,j+1),(\alpha-1,\beta+1)$
satisfy the same set of conditions out of 1--7 (with inequalities relaxed).
\qed
\enddemo

Denote $n-k$  by  $m$. If we arrange variables  $\ln t_{ij}$  into a vector
$$
{\tilde t} =(\ln t_{11},\ldots,\ln t_{1 m},\ldots, \ln t_{k1},\ldots,\ln t_{k m}),
\tag 3.19
$$
the coefficient
matrix $\Omega_{km}$ of Poisson brackets (3.18) will look as follows:
$$
\Omega_{km}=\pmatrix
A_m & \one_m & \one_m & \cdots & \one_m\\
-\one_m & A_m & \one_m & \cdots & \one_m\\
-\one_m & -\one_m & A_m & \cdots & \one_m\\
\vdots & \vdots &\vdots &\ddots &\vdots \\
-\one_m & -\one_m & -\one_m & \cdots& A_m
\endpmatrix
 = A_m\otimes \one_k  - \one_m\otimes A_k,
\tag 3.20
$$
where $A_1=0$ and
$$
A_{m}=\Omega_{1m}^T=
\pmatrix
0 & -1 & -1 & \cdots & -1\\
1 & 0& -1 & \cdots & -1 \\
1 & 1  & 0 & \cdots & -1 \\
\vdots & \vdots &\vdots &\ddots &\vdots \\
1 & 1 & 1 &\cdots & 0
\endpmatrix.
$$

We now proceed to compute a maximal dimension
of a symplectic leaf of the bracket (3.11).

First note that left  multiplying $(\lambda \one_m  + A_m^T)$  by
$C_m= \one_m + e_{1m} - \sum_{i=2}^m e_{i,i-1}$, where $e_{ij}$ is
a $(0,1)$-matrix with a unique $1$ at position $(i,j)$,
results in the following matrix:
$$
B_{m}(\lambda)= \pmatrix
\lambda -1 & 0 & 0& \cdots & 0 & \lambda +1\\
-\lambda -1 & \lambda -1& 0 & \cdots & 0 &0\\
0 & -\lambda -1 & \lambda -1& \cdots & 0 &0 \\
\vdots & \vdots &\vdots &\ddots &\vdots &\vdots\\
0 &0 &0 &\cdots & -\lambda -1 & \lambda -1
\endpmatrix.
$$
Since the determinant of $B_{m}(\lambda)$ is easily computed
to be equal to $(\lambda -1)^m + (\lambda +1)^m$, it follows that
the spectrum of $A_m$ is given by
$$
\left (\frac{\lambda +1}{\lambda -1}\right )^m = -1.
\tag 3.21
$$

Next, we observe that using block row transformations similar to row
transformations
applied to $(\lambda \one_m  + A_m^T)$ above, one can reduce $\Omega_{km}$
to a matrix
$B_k(A_m)$ obtained from $B_k(\lambda)$ by replacing $\lambda$ with $A_m$ and
$1$ with $\one_m$. Since  $(A_m - \one_m)$ is invertible by (3.21),
we can left multiply $B_k(A_m)$ by diag$((A_m - \one_m)^{-1},\ldots, (A_m - \one_m)^{-1})$
and conclude that the kernel of $\Omega_{km}$ coincides with that of the matrix
$$
\pmatrix
\one_m & 0 &0 & \cdots & W\\
-W & \one_m& 0 & \cdots & 0 \\
0 & -W & \one_m& \cdots & 0 \\
\vdots & \vdots &\vdots &\ddots &\vdots \\
0 & 0 & 0& \cdots & \one_m
\endpmatrix,
$$
where $W=(A_m+\one_m)(A_m-\one_m)^{-1}$.
The kernel consists of vectors of the form
$\pmatrix v & Wv & \ldots &W^{k-1}v\endpmatrix^T$, where $v$ satisfies
condition $(W^k + \one_m)v=0$. In other words,
the kernel of $\Omega_{km}$ is parametrized by the $(-1)$-eigenspace
of $W^k$. Due to (3.21) the dimension of this eigenspace
is equal to
$$
d_{km}= \# \left \{\nu \in  \Bbb{C}\: \nu^k=\nu^m = -1 \right \}.
$$
Moreover, it is not hard to check that
$$
W=(A_m+\one_m)(A_m-\one_m)^{-1} = -e_{m1} + \sum_{i=2}^m e_{i-1,i},
\tag 3.22
$$
and therefore $W^k=\sum_{i=k+1}^m e_{i-k,i} -\sum_{i=1}^{k} e_{i+m-k,i}$.
A $(-1)$-eigenspace of $W^k$ is non-trivial if and only if there exist natural
numbers
$p,q$ such that $(2p-1)k=(2q-1)m$ (we can assume
$(2p-1)$ and $(2q-1)$ are co-prime). Let $l=\gcd(k,m)$,
then every non-trivial $(-1)$-eigenvector of $W^k$
is a linear combination of vectors $v(i)=(v(i)_j)_{j=1}^m$, $i\in [1,l]$,
that can be described as follows: $v(i)_{i+\alpha l} =(-1)^{\alpha}$ for
$\alpha=0,\ldots, \frac{m}{l} -1$, and all the other entries of $v(i)$ vanish.
To analyze the corresponding element
$\pmatrix v(i)& Wv(i) &\ldots & W^{k-1}v(i)\endpmatrix^T$
of the kernel of $P_{km}$, we represent it as a  $k\times m$ matrix
$$
V(i)=\pmatrix v(i)\\ Wv(i)\\ \vdots\\ W^{k-1}v(i)\endpmatrix.
$$

From the form of $W$ one concludes that $V(i)$ is a matrix
of $0$s and $\pm 1$s that has a Hankel structure, i.e. its entries do not
change along anti-diagonals. More precisely,
$$
V(i)_{pq} =\left\{\alignedat2
(-&1)^{\alpha} \qquad  &\text{if $p+q =i + \alpha l$}, \\
&0 \qquad  &\text{otherwise}\endalignedat
\right.
$$
Here $\alpha$ changes from $0$ to $\frac{n}{l}-1 $.
Since to each element $V$ of the kernel of $\Omega_{km}$ (represented
as a $k\times m$ matrix) there corresponds a Casimir function $I_V$ of
$\omega_{G(k,n)}$
given by $I_V =\prod_{i=1,j=1}^{k,m} t_{ij}^{V_{ij}}$, the observation above
together with (3.17) implies that  {\it on an open dense subset of $G(k,n)$
the algebra of Casimir functions is generated by monomials in}
$$
J_1=F_{11},\ldots,J_k=F_{k1}, J_{k+1}=F_{k2},\ldots,J_n=F_{km}.
\tag 3.23
$$
In particular,
$$
I_{V(i)} = \prod_{\alpha=0}^{\frac{n}{l}-1} J_{i+\alpha l}^{(-1)^\alpha}.
\tag 3.24
$$

Thus we have proved

\proclaim{Theorem 3.5} Let $l=\gcd(k,n)$. The co-dimension of a maximal
symplectic leaf of $G(k,n)$ is equal to $0$ if $\frac{k}{l}$ is even
or $\frac{n-k}{l}$ is even, and is equal to $l$ otherwise.
In the latter case, a symplectic leaf via a point in general position
is parametrized by values of Casimir functions $I_{V(i)}$, $i\in [1,l]$,
defined in {\rm (3.24)}.
\endproclaim

\subheading{3.3. Cluster algebra structure on Grassmannians compatible with the
Poisson bracket}
Our next goal is to build a cluster algebra $\AA_{G(k,n)}$
associated with the Poisson bracket
(3.11). The initial cluster consists of functions
$$\multline
  f_{ij} =(-1)^{(k-i)(l(i,j)- 1)} F_{ij}=
\frac{\pi_{([1,k]\setminus [i-l(i,j),\ i])
\cup [ j+k, j+l(i,j)+k]}}{\pi_{[1,k]}},\\
i\in [1,k],\quad j\in [1, n-k]).
\endmultline\tag 3.25
$$
We designate functions $f_{11},f_{21},\ldots, f_{k1}, f_{k2},\ldots, f_{km}$
(cf. (3.23))  to serve as tropical coordinates. This choice is motivated by
the last statement
of Theorem~3.5 and by the following observation: let $I,J$ be the row and
column sets of the minor that represents one of the functions (3.23), then
 for any pair $(\alpha, \beta)$,
$\alpha\in [1, k]$, $\beta\in [1,m]$, condition (3.13) is satisfied
and thus, functions (3.23) have log-canonical brackets with all coordinate
functions $y_{\alpha \beta}$.

Now we need to define the matrix $Z$ that gives rise to cluster
transformations compatible with the Poisson structure. We want to choose $Z$
in such a way that the submatrix of $Z$ corresponding to cluster coordinates
will be skew-symmetric and irreducible. According to (1.7) and to our choice
of tropic
coordinates, this means that $Z$ must satisfy
$$
 Z \Omega^F =  const \cdot  \pmatrix \diag(P,\ldots, P) & 0\endpmatrix
$$
where $\P= \sum_{i=1}^{m-1} e_{i,i+1}$ is a $(m-1)\times m$ matrix
and $\Omega^F$ is the coefficient matrix of Poisson brackets $\omega$
in the basis $F_{ij}$.

Let ${\tilde t}$ be defined as in (3.19), and let
$$
\tilde F=(\ln F_{11},\ldots,\ln F_{1 m},\ldots,\ln F_{k1},\ldots,\ln F_{k m}).
$$
Then $ {\tilde t}=J \tilde F^T$, where
$$
J=
\pmatrix
  \one_m & 0 & \cdots & 0 & 0\\
  -S & \one_m & \cdots & 0 &0 \\
   \vdots & \vdots & \ddots & \vdots &\vdots \\
  0 &0 &\cdots & -S & \one_m
\endpmatrix
$$
and $S= \sum_{i=2}^m e_{i-1,i}$.
Then $\Omega^F = J \Omega_{km} J^T$.

Define a $(k-1)m\times km$ block bidiagonal matrix
$$
V=
\pmatrix
  \one_m - W^{-1} & W^{-1}-\one_m & \cdots & 0 &0 \\
  0 & \one_m - W^{-1} &  \cdots & 0 &0\\
  \vdots & \vdots & \ddots & \vdots & \vdots  \\
  0 & 0 & \cdots &  \one_m - W^{-1} & W^{-1}-\one_m
\endpmatrix.
$$
Observe that $P$ is the upper $(m-1)\times m$ submatrix of $S$, and
$P W^{-1}= P S^T$.

Since by (3.20), (3.22),
$$
V \Omega_{km}= 2 \
\pmatrix
  \one_m & - W^{-1}& 0 & \dots &0 & 0 \\
  0 & \one_m & - W^{-1}& \dots & 0 &0\\
  \vdots & \vdots & \vdots & \ddots & \vdots &\vdots \\
  0 &0& 0&\dots &   \one_m  & - W^{-1}
\endpmatrix,
$$
we obtain
$$
\frac{1}{2}\diag(P,\ldots, P) V\Omega_{km}=
\pmatrix\diag(P,\ldots, P) &0\endpmatrix J^T.
$$
Define
$$
Z = \diag(P,\ldots, P) V J=
\pmatrix
  Z_0 & Z_{1}& 0 & \dots & 0 \\
  Z_{-1} & Z_0 & Z_1& \dots & 0 \\
  0 & Z_{-1} & Z_0 & \dots & 0 \\
  \vdots & \vdots & \vdots & \ddots & \vdots  \\
  0 & 0 &0 &\dots &   Z_0
\endpmatrix,
$$
where $Z_0=P(\one_m-W^{-1})(\one_m+S)$, $Z_1=-P(\one_m-W^{-1})$, and
$Z_{-1}==-P(\one_m-W^{-1})S$.
Then
$$
Z \Omega^F =  \diag(P,\ldots,P) V\Omega_{km} (J^T)^{-1}=
2\pmatrix\diag(P,\ldots, P) &0\endpmatrix,
$$
as needed. Note  that for
$x =(x_{11},\ldots,x_{1,m},\ldots, x_{k1},\ldots,x_{k,m})$ one has
$$
(Z x)_{ij} =
x_{i+1,j} + x_{i, j-1} + x_{i-1,j+1} - x_{i+1,j-1} - x_{i, j+1} - x_{i-1,j}.
$$
It is easy to see that the submatrix of $Z$ corresponding to the non-tropic
coordinates is indeed skew-symmetric and irreducible.

The matrix $Z$ thus obtained can be conveniently represented by a
directed graph with vertices
forming a rectangular $k\times (n-k)$ array and labeled by pairs of integers
$(i,j)$, and edges $(i,j) \to (i,j+1),\ (i+1,j) \to (i,j)$ and $(i,j) \to (i+1,j-1)$
(cf.~Fig.~1).

   \vskip 15pt
\centerline{\hbox{\epsfysize=3.5cm\epsfbox{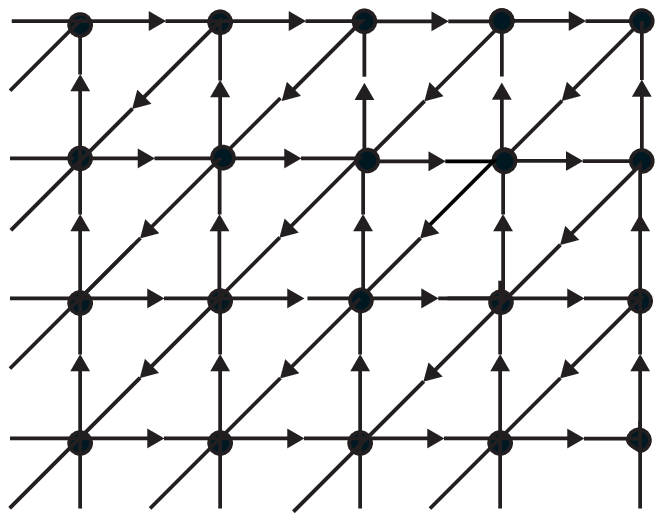}}}
\midspace{0.1mm} \caption{Fig.~1. Graph that corresponds to $G(k,n)$}


\proclaim{Proposition 3.6}
For every $i\in [1,k-1]$ and every $j\in [2,n-k]$, a cluster
variable $f'_{ij}$ obtained via the cluster transformation from the initial
cluster {\rm (3.25)} is a regular function on $G_0(k,n)$.
\endproclaim

\demo{Proof} The proof is based on Jacobi's determinantal identity
$$
A_{(\alpha\beta I)}^{(\gamma\delta J)} A_I^J=
A_{(\alpha I)}^{(\gamma J)}A_{(\beta I)}^{(\delta J)}
-A_{(\alpha I)}^{(\delta J)}A_{(\beta I)}^{(\gamma J)}\ .
\tag 3.26
$$

We will consider the following cases:

(i) $ 1< i < n-k+1-j $. In this case
$f_{ij} = (-1)^{(k-i)i} Y_{[1,i]}^{[j,j+i-1]}$, and
(3.26), (3.25) imply
$$\multline
(-1)^{i(k+i)} f_{i+1,j-1}  f_{i,j+1} f_{i-1,j}= Y_{[1,i+1]}^{[j-1,j+i-1]} Y_{[1,i-1]}^{[j,j+i-2]}
Y_{[1,i]}^{[j+1,j+i]}\\
= \left (  Y_{[1,i-1]\cup i+1}^{[j,j+i-1]} Y_{[1,i]}^{[j-1,j+i-2]}-
Y_{[1,i-1]\cup i+1}^{[j-1,j+i-2]} Y_{[1,i]}^{[j,j+i-1]}  \right )Y_{[1,i]}^{[j+1,j+i]}
\endmultline
$$
and
$$\multline
(-1)^{i(k+i)}f_{i+1,j}  f_{i,j-1} f_{i-1,j+1}= Y_{[1,i+1]}^{[j,j+i]} Y_{[1,i-1]}^{[j+1,j+i-1]}
Y_{[1,i]}^{[j-1,j+i-2]}\\
= \left (  Y_{[1,i-1]\cup i+1}^{[j+1,j+i]} Y_{[1,i]}^{[j,j+i-1]}-
Y_{[1,i-1]\cup i+1}^{[j,j+i-1]} Y_{[1,i]}^{[j+1,j+i]}\right )Y_{[1,i]}^{[j+1,j+i]}.
\endmultline
$$
Therefore,
$$\multline
f'_{ij} =
\frac{f_{i+1,j-1}f_{i,j+1}f_{i-1,j} +f_{i+1,j}f_{i,j-1}f_{i-1,j+1}}{f_{ij}}\\
=  Y_{[1,i-1]\cup i+1}^{[j+1,j+i]} Y_{[1,i]}^{[j-1,j+i-2]}-
Y_{[1,i-1]\cup i+1}^{[j-1,j+i-2]} Y_{[1,i]}^{[j+1,j+i]}.
\endmultline
$$

Other cases can be treated similarly. Below, we present corresponding
expressions for $f'_{ij}$.

(ii) $ n-k+1-i < j < n-k$. Then
$$
f'_{ij} =
Y^{j-1\cup [j+1,m]}_{[\alpha+1,i+1]} Y_{[\alpha-1,i-1]}^{[j,m]}-
Y^{j-1\cup [j+1,m]}_{[\alpha-1,i-1]} Y_{[\alpha+1,i+1]}^{[j,m]},
$$
where $m=n-k$ and $\alpha=i=j=k-n$.

(iii) $ n-k+1-i = j < n-k$. Then
$$
f'_{ij} =(-1)^{k-i}
\left ( Y^{[j-1,m]}_{[1,i+1]} Y_{[2,i-1]}^{[j+1,m-1]}-
Y^{[j-1,m-1]}_{[2,i+1]} Y_{[1,i-1]}^{[j+1,m]}\right ).
$$

(iv) $ i = 1$, $j < n-k$. Then $f'_{1j}= Y_{12}^{j-1,j+1}$.

(v) $ i > 1$, $j =n-k$. Then $f'_{1j}= - Y_{i-1,i+1}^{n-k-1,n-k}$.

(vi) $ i = 1$, $j = n-k$. Then  $f'_{1,n-k}=(-1)^{k} y_{2, n-k-1}$.

In all six cases, $f'_{ij}$ is a polynomial in variables $y_{pq}$, which
proves the assertion.
\qed
\enddemo

\proclaim{Proposition 3.7}
For every $i\in [1,k-1]$ and every $j\in [2,n-k]$,
the coordinate function
$y_{ij}$ belongs to some cluster obtained from the initial one.
\endproclaim

\demo{Proof} Let $m=n-k$ as before, and let $\bar Y$ be the
matrix obtained from $Y$ by deleting
the first row and the last column. Denote by $\bar F_{ij}$, $\bar f_{ij}$,
$i\in [2,k]$, $j\in [1,n-k-1]$, the functions defined by (3.16),
(3.25) with
$Y$ replaced by  $\bar Y$. Define also $f = (f_{ij})$ and  $\bar f =
(\bar f_{ij})$.

Let us consider the following composition of cluster transformations:
$$
T=T_{k-1}\circ \cdots \circ T_{1},
\tag 3.27
$$
where
$$
T_\gamma= T_{k-1, m-\gamma+1} \circ \cdots \circ T_{\gamma+1, m-\gamma+1}  \circ T_{\gamma\ 2} \circ \cdots \circ T_{\gamma, m-\gamma+1}
$$
for $\gamma=1,\ldots, k-1$.
Note that every cluster transformation $T_{ij}$, $i=2,\ldots,k$,
$j=1,\ldots,n-k-1 $, features in (3.27) exactly once.

We claim that
$$\alignedat2
(Tf)_{ij}&=\bar f_{ij},&\quad &i\in [1,k-1]; j=2,m-1, \\
(TZ)_{(ij),(\alpha \beta)}&= Z_{(ij),(\alpha \beta)},
&\quad &i,\alpha\in [1,k-1]; j,\beta=2,m-1; j+\beta > 2; i+\alpha < 2k.
\endalignedat
\tag3.28
$$
In particular, $(Tf)_{1 j} = Y_{2,k-1}$, $j\in [2,m]$,
and $(Tf)_{i m} = Y_{i+1,m-1}$, $j\in [2,m]$. Applying the same
strategy to $\bar Y$ etc., we will eventually recover
all matrix entries of $Y$.

To prove (3.28), it is convenient to work with graphs associated with
the matrix $Z$ and its transformations, rather
that with matrices themselves. Using the initial graph given on Fig.~1
and using the remark at the end of \S~1.1, it is not hard to convince
oneself that
at the moment when $T_{ij}$ (considered as a part of the composition
(3.27))
is applied, the corresponding graph changes according to Fig.~2.
In the latter figure, the white circle denotes the  vertex $(i,j)$ and only
the vertices
connected with  $(i,j)$ are shown. If $i=1$ (resp., $j=m$) then vertices above
(resp. to the right of)  $(i,j)$ and edges that connect them to $(i,j)$
should be ignored.

   \vskip 15pt
\centerline{\hbox{\epsfxsize=5cm\epsfbox{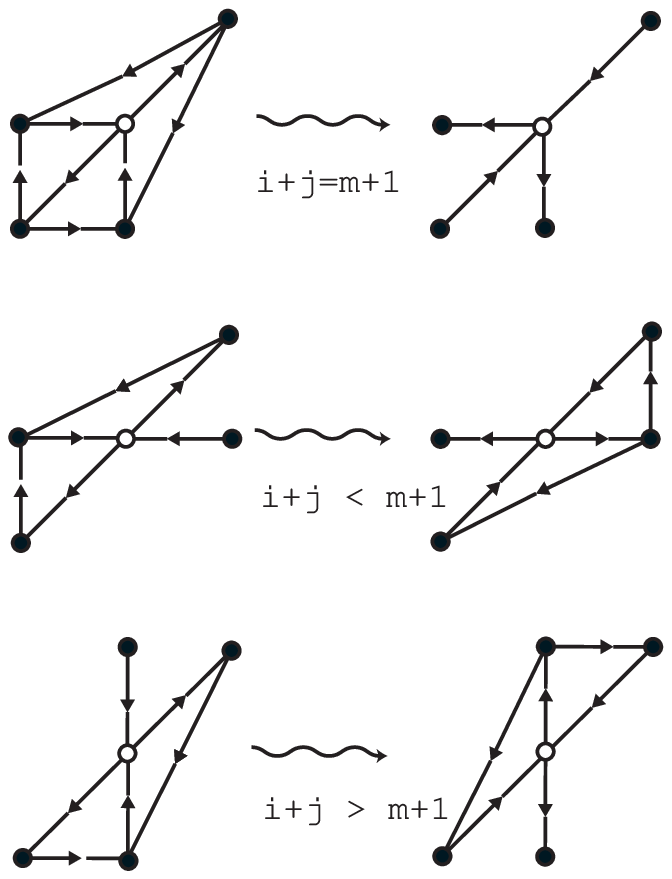}}}
\midspace{0.1mm} \caption{Fig.~2. Transformation $T_{ij}$}
\vskip 15pt

It follows from Fig.~2 that

(i) the direction of an edge between
$(\alpha,\beta)$ and $(\alpha+1,\beta-1)$ changes when $T_{\alpha \beta}$ is applied and is restored when $T_{\alpha+1,\beta-1}$ is applied;

(ii) a horizontal
edge between $(\alpha,\beta)$ and $(\alpha,\beta+1)$ is erased when
$T_{\alpha-1,\beta-1}$ is applied
and is restored with the original direction when  $T_{\alpha+1,\beta}$
is applied;

(iii) a vertical
edge between $(\alpha,\beta)$ and $(\alpha+1,\beta)$ is erased when
$T_{\alpha,\beta+1}$ is applied
and is restored with the original direction when  $T_{\alpha-1,\beta+1}$
is applied;

(iv) an edge between $(\alpha,\beta+1)$ and $(\alpha+1,\beta-1)$ is
introduced when  $T_{\alpha \beta}$ is applied
and is erased when  $T_{\alpha+1,\beta}$ is applied;

(v) an edge between $(\alpha-1,\beta)$ and $(\alpha+1,\beta-1)$ is
introduced when  $T_{\alpha \beta}$ is applied
and is erased when  $T_{\alpha,\beta-1}$
is applied.

Thus, after applying all transformations that constitute $T$ in (3.27),
one obtains a directed graph whose upper right $(k-1)\times (m-1)$ part
(not taking into account edges between the vertices in the first column and the
last row) coincides with that of the initial graph on Fig.~1.
This proves the second equality in (3.28).

To prove the first equality in (3.28), note that, by the definition
of (3.27), a cluster coordinate $f_{ij}$ changes to $(Tf)_{ij}$
at the moment when $T_{ij}$ is applied and stays unchanged afterwards.
In particular, on Fig.~2, coordinates that correspond to vertices above
and to the right
of $(i,j)$ are coordinates of $Tf$, while
coordinates that correspond to vertices above and to the right
of $(i,j)$ are coordinates of $f$.
Thus,  (3.28) will follow if we show that
$$\align
f_{ij} \bar f_{ij}&= f_{i+1,j-1} \bar f_{i-1,j+1} + f_{i,j-1}\bar f_{i+1,j},
\quad i+j=m+1,\\
f_{ij} \bar f_{ij}&= f_{i+1,j-1} \bar f_{i-1,j+1} + f_{i,j-1} \bar f_{i,j+1},
\quad i+j < m+1,\\
f_{ij} \bar f_{ij}&= f_{i+1,j-1} \bar f_{i-1,j+1} + f_{i+1,j} \bar f_{i-1,j},
\quad i+j > m+1.
\endalign
$$
But referring to definitions of $f$ and $\bar f$, one finds that the three
equations
above are just another instances of Jacobi's identity
(3.26). The proof is complete.
\qed
\enddemo

\subheading{3.4. The number of connected components of refined Bruhat
cells in real Grassmannians}
Consider the union of regular $\R^*$-orbits in $\X^0_{G(k,n)}$
corresponding to the described above
cluster algebra $\AA_{G(k,n)}$ compatible with the Sklyanin Poisson
bracket in  $G(k,n)$.
Recall that  by construction functions
$f_{11},\dots,f_{k1},f_{k2},\dots,f_{km}$ serve as tropical coordinates.
Moreover, any matrix element in the standard representation of the maximal
Bruhat cell in the Grassmannian enters as a cluster coordinate for some
cluster in $\AA_{G(k,n)}$. Therefore,
$\X^0_{G(k,n)}$ is naturally embedded into $G(k,n)$ and we can consider it
as the subset in $G(k,n)$ determined by the conditions that all
tropical coordinates $f_{ik}$ and $f_{kl}$ do not vanish.
Tropical coordinates are all ``cyclically dense'' minors among
all the Pl\"ucker coordinates, i.e., minors containing all columns with
indices $i,i+1,\dots,i+k$ or
$i,i+1,\dots,i+l=n,1,2,\dots, k+i+1-n$.
We call $\X^0_{G(k,n)}$ a {\it refined open Bruhat cell\/} in $G(k,n)$, since it is
an intersection of $n$ usual open Bruhat cells in $G(k,n)$ in general position.

A method to compute the number of connected components of $\X^0_{G(k,n)}$ was
discussed in Section~2 of this paper.
Let us recall certain notions and results from~\cite{SSVZ}.

We denote by $\gr$ the graph corresponding to the matrix $Z$ (see the
remark at the end of Section~1.1), by $\F_2^\gr$ the vector space
generated by the characteristic vectors of the vertices of $\gr$,
and by $\eta_{\gr}$ the corresponding skew-symmetric bilinear form on
$\F_2^\gr$ (in our case, $\eta_{\gr}(e_i,e_j)=z_{ij}$). Similarly,
$\F_2^C$ denotes a subspace of $\F_2^\gr$ generated by the vertices
corresponding to cluster variables.
A finite (undirected) graph is said to be {\it $E_6$-compatible\/} if it is
connected and contains an induced subgraph with $6$ vertices isomorphic to
the Dynkin graph $E_6$. A directed graph is said to be $E_6$-compatible if
the corresponding undirected graph obtained by replacing each directed edge
by an undirected one is  $E_6$-compatible.

\proclaim{Theorem 3.8 {\rm (\cite{SSVZ, Th.~3.11})}} Suppose that the
induced subgraph of $\gr$ spanned by the vertices corresponding
to cluster variables
is $E_6$-compatible. Then the number of $\Gamma$-orbits in
$\F_2^\gr$ is equal to
$$
2^t\cdot(2+2^{\dim(\F_2^C\cap \ker \eta_{\gr})}),
$$
where $t$ is the number of tropic variables.
\endproclaim

Combining this theorem with Theorem 2.10 we get the following corollary.

\proclaim{Corollary 3.9} The number of
connected components of a refined open Bruhat cell $\X^0_{G(k,n)}$ is equal to
$3\cdot 2^{n-1}$ if $k>3$, $n-k>3$.
\endproclaim

\demo{Proof}
Indeed, by Theorem~2.10 we know that the number of connected components of
$\X^0_{G(k,n)}$ equals the number
of orbits of $\Gamma$-orbits in $\F_2^\gr$, where the graph $\gr$ is shown on Fig.~1,
and the subset $C$ is formed by all the vertices except for the first column
and the last row.
Since in the case $k\gs 4$, $n\gs 8$ the subgraph spanned by $C$
is evidently $E_6$-compatible, Theorem~3.8 implies
that to prove the statement it is enough to show that
$\F_2^C\cap \ker \eta_{\gr}=0$;
in other words, that there is no nontrivial vector in  $\ker \eta_{\gr}$
with vanishing
tropical components.

Indeed, let us denote such a vector by $h\in \F_2^C$, and let $\delta_{ij}$
be the $i,j$-th basis vector
in $\F_2^\gr$.
Note that the condition $ \eta_{\gr}(h, \delta_{k,n-k})=0$ implies that
$h_{k-1,n-k}=0$.
Further, assuming $h_{k-1,n-k}=0$ we see that the condition
$ \eta_{\gr}(h, \delta_{k,n-k-1})=0$
implies $h_{k-1,n-k-1}=0$ and so on.
Since $\eta_{\gr}(h, \delta_{kj})=0$ for any $j\in[1,m]$,
we conclude that $h_{k-1,j}=0$ for any $j\in[1,n-k]$. Proceeding by induction we prove
that $h_{ij}=0$ for any $i\in[1,k]$ and any $j\in[1,n-k]$. Hence $h=0$.
Note that $t=n-1$, which accomplishes the proof of the statement.
\qed
\enddemo

It is easy to notice that the number of connected components for
$\X^0_{G(k,n)}$
equals the number of connected components for  $\X^0_{G(n-k,n)}$.
Therefore, taking in account Corollary~3.9, in order to find the number of connected components
for refined open Bruhat cells for all Grassmannians we need to consider
only two remaining cases: $G(2,n)$ and $G(3,n)$.

\proclaim{Proposition 3.10} The number of connected components of a refined
open Bruhat cell equals to $(n-1)\cdot 2^{n-2}$ for $G(2,n)$, $n\gs3$, and to
$3\cdot 2^{n-1}$ for $G(3,n)$, $n\gs6$.
\endproclaim

\demo{Proof} The proof  follows immediately from Lemmas 1,2 and corollary of
Lemma~3 of~\cite{GSV}.
Following notations of~\cite{GSV}, let us denote by $U$ the subgraph of the
graph $\gr$ (for the Grassmannian $G(k,n)$)
consisting of  the first $k-1$ rows,  and by $L$ the subgraph containing
only the last row.
The corresponding vector subspaces
of $\F_2^\gr$
are denoted by $\F_2^U$ and $\F_2^L$; the corresponding
subgroups of $\Gamma$ are denoted by $\Gamma^U$ and $\Gamma^L$.
Lemma~1 of~\cite{GSV} states that for any $G(k,n)$ we have
$2^{n-k}\#(\Gamma^U)$ orbits of
$\Gamma$-action, where $2^{n-k}$ is the number of vectors in
$\F_2^L$, and $\#(\Gamma^U)$
is the number of $\Gamma^U$-orbits
in $\F_2^U$.
Lemma~2 counts   $\#(\Gamma^U)=n-1$ for  $G(2,n)$, $n\gs3$.
Finally, the corollary of Lemma~3 calculates  $\#(\Gamma^U)=12$ for  $G(3,n)$,
$n\gs6$.
Substituting these values of $\#(\Gamma^U)$  into the above mentioned
formula from Lemma~1  proves the proposition.\qed
\enddemo

\end